\documentclass[11pt]{article}

\usepackage{amssymb,latexsym,amsmath,mathrsfs}
\usepackage{graphicx}
\usepackage{epstopdf}
\usepackage{hyperref}
\usepackage{booktabs}
\usepackage[table]{xcolor}
\numberwithin{equation}{section}

\newcommand{\be}{\begin{equation}}
  \newcommand{\ee}{\end{equation}}

\newtheorem{theorem}{Theorem}[section]

\newtheorem{lem}{Lemma}[section]

\newcommand{\att}[1]{{
		{#1}}} 

\begin{document}

\title{A Laguerre homotopy method for optimal control of nonlinear systems in semi-infinite interval}

\author{ Haijun Yu  $^{1,2}$, Hassan Saberi Nik $^\dag$ \\
$^{1}$ \small{School of Mathematical Sciences, University of Chinese Academy of Sciences, Beijing 100049, China}\\
$^{2}$\small{NCMIS \& LSEC, Institute of Computational Mathematics and Scientific/Engineering Computing, }\\
\small{Academy of Mathematics and Systems Science, Beijing 100190, China}\\
\small{email: hyu@lsec.cc.ac.cn}\\
$^\dag$\small{Department of Mathematics, Neyshabur Branch, Islamic Azad University, Neyshabur, Iran}\\
\small{email: saberi\underline{~}hssn@yahoo.com}
\\}

\maketitle

\begin{center}
  \large{\bf Abstract}
\end{center}

This paper presents a Laguerre homotopy method for optimal control problems in semi-infinite intervals (LaHOC), with particular interests given to nonlinear interconnected large-scale dynamic systems. In LaHOC, spectral homotopy analysis method is used to derive an iterative solver for the nonlinear two-point boundary value problem derived from Pontryagins maximum principle. A proof of local convergence of the LaHOC is provided. Numerical comparisons are made between the LaHOC, Matlab BVP5C generated results and results from literature for two nonlinear optimal control problems. The results show that LaHOC is superior in both accuracy and efficiency.\\


\section{Introduction}
Large-scale systems are found in many practical applications, such as power systems and physical plants. During the  past several years, the problem of analysis and synthesis for dynamic large-scale systems has received considerable attention. Based on the characteristics of large-scale systems many results have been proposed, such as modelling, stability, robust control, decentralized, and so on \cite{Jamshidi,Wasim,Holland,Huang,Chen,Yan}.

The optimal control problem (OCP) of nonlinear large-scale systems has been widely investigated in recent decades.  For instance, a new successive approximation approach
(SAA) was proposed in \cite{Tang1}. In this approach, instead of directly solving the nonlinear large-scale two-point boundary value problem (TPBVP), derived from the
maximum principle, a sequence of non-homogeneous linear time-varying TPBVPs is solved iteratively.  Also, in \cite{Jajarmi} a new technique, called the
modal series method, has been has been extended to solve a class of infinite horizon OCPs of nonlinear interconnected large-scale dynamic systems,
where the cost function is assumed to be quadratic and decoupled. This method provides the solution of autonomous nonlinear systems
in terms of fundamental and interacting modes. Conventional methods of optimal control are generally impractical for many nonlinear large-scale systems because of the dimensionality problem and high complexity in calculations. One example is the state-dependent Riccati equation (SDRE) method \cite{Chang}. Although this scheme has been widely used in many applications, its major limitation is that it needs to solve a sequence of matrix Riccati algebraic equations at each sample state along the trajectory. This property may take a long computing time and large memory space. Therefore, developing new methods is necessary for solving nonlinear large-scale optimal control problems \cite{Kaffash}.

The use of spectral methods for optimal control problem usually leads to more efficient method than finite element or finite different approaches. Chebyshev and Legendre method are commonly used for problems in finite intervals \cite{RossF2003, RossF2004}. For infinite or semi-infinite intervals, there are several choices for the approximation bases: Hermite polynomials/functions \cite{Tang93}, Laguerre polynomials /functions \cite{GuoWangTianWang2008}, mapped Jacobi bases \cite{ShenW2009, ShenY2012, ShenWY2014}. Furthormore, one class of very important applications of OCP in unbounded intervals is the minimum action method (MAM) \cite{ERV2004} used in finding the most probable transition path in phase transition phenomena. Using MAM to study spatial extended transitions, such as
fluid instability transition is usually equivalent to solve a large-scaled nonlinear optimal control problem \cite{WanYE2015, WanY2017}.

The homotopy analysis method is an analytical technique for solving nonlinear differential equations. The HAM \cite{Liao3,Liao4} was first proposed by Liao in 1992 to solve lots of nonlinear problems. This method has been successfully applied to many nonlinear problems, such as physical models with an infinite number of singularities \cite{Liao5},  nonlinear eigenvalue problems \cite{Abbasbandy1}, fractional Sturm-Liouville problems \cite{Abbasbandy2}, optimal control problems \cite{Saberi3,Jajarmi2},  Cahn-Hilliard initial value problem \cite{Vangorder1}, semi-linear elliptic boundary value problems \cite{Vangorder2} and so on \cite{Zhu}. 
The HAM contains a certain auxiliary parameter $\hbar$ which provides us with a simple way to adjust and control the convergence region and rate of convergence of the series solution. Moreover, by means of the so-called $\hbar$-curve, it is easy to determine the valid regions of $\hbar$ to gain a convergent series solution.
The HAM however suffers from a number of restrictive measures, such as the requirement that the solution sought ought to conform to the so-called rule of solution expression and the rule of coefficient ergodicity. These HAM requirements are meant to ensure that the implementation of the method results in a series of differential equations which can be solved analytically.

Recently, Motsa et al. \cite{Motsa1, Motsa2,Saberinik} proposed a spectral modification of the homotopy analysis method, the spectral-homotopy analysis method (SHAM). The SHAM approach imports some of the ideas of the HAM such as the use of the convergence controlling auxiliary parameter. In the implementation of the SHAM, the sequence of the so-called deformation differential equations are converted into a matrix system by applying Chebyshev or Legendre pseudospectral method \cite{Saberinik}.  But so far, to our knowledge, there is no work concerning the combination of Laguerre polynomails \cite{JieShen} with the HAM. This paper presents a spectral homotopy analysis method based on modified Laguerre-Radau interpolation to solve nonlinear large-scale optimal control problems. This process has several advantages. First, it possesses the spectral accuracy \cite{Saadatmandi,Bhrawy}. Next, it is easier to be implemented, especially for nonlinear systems. Furthermore, it is applicable to long-time calculations.

The paper is organized as follows. The nonlinear interconnected OCP and optimality conditions is described in section 2. In Section 3, we propose the new algorithm by using the modified Laguerre polynomials.  The convergence of the proposed method is proved in section 4. We present numerical results in Section 5, which demonstrate the spectral accuracy of proposed methods. The final section is for concluding remarks.

\section{The nonlinear interconnected OCP}
Consider a nonlinear interconnected large-scale dynamic system which can be decomposed into $N$ interconnected subsystems. The $ith$ subsystem for $i =1,2,\cdots, N$ is
described by:
\begin{equation}
  \label{Eq:Nonlinear dynamical system}
  \begin{array}{l}
    \dot x_i(t)=A_ix_i(t)+B_iu_i(t)+f_i(x(t)),\quad  t> t_0, \\
    x_i(t_0)=x_{i_0},
  \end{array}
\end{equation}
with $x_i\in \mathbb{R}^{n_i}$ denoting the state vector, $u_i\in\ \mathbb{R}^{m_i}$ the control vector of the ith subsystem, respectively, $x=(x_1^T,x_2^T,\cdots x_N^T)^T,$ $\displaystyle\sum_{i=1}^{N}n_i=n,$ $F_i:\mathbb{R}^{n}\rightarrow \mathbb{R}^{n_i}$ is a nonlinear analytic vector function where $F_i (0)=0,$ and
$x_{i_0} \in \mathbb{R}^{n_i}$ is the initial state vector. Also, $A_i$ and $B_i$ are constant matrices of appropriate dimensions such that the pair $(A_i , B_i )$ is
completely controllable \cite{Jajarmi}. Furthermore, the infinite horizon quadratic cost function to be minimized is given by:

\begin{equation}
  \label{Eq:Performance index}
  J=\frac{1}{2}\sum_{i=1}^{K}\left\{\int_{t_0}^{\infty} (x_i^T(t)Qx_i(t)+u_i^T(t)R_iu_i(t))dt\right\}
\end{equation}
\noindent where $Q_i \in \mathbb{R}^{n_i\times n_i}$ and $R_i\in \mathbb{R}^{m_i\times m_i}$ are positive semidefinite and positive definite matrices, respectively.
Note that the quadratic cost function (\ref{Eq:Performance   index}) is assumed to be decoupled as a superposition of the cost functions of the subsystems.

According to Pontryagin's maximum principle, the optimality conditions are obtained as the following nonlinear TPBVP:
\begin{equation}
  \label{Eq:Optimality conditions 2}
  \begin{array}{l}
    \dot x_i(t) =A_ix_i(t)-B_iR_i^{-1}B_i^T\lambda_i(t)+f_i(x(t)), \quad t> t_0,\\
    \dot \lambda_i(t) = -Q_ix_i(t)-A_i^T\lambda_i(t)-\Psi_i(x(t), \lambda(t)), \quad t> t_0,\\
    x_i(t_0)=x_{i_0}, \lambda_i(\infty)=0,\\
    i=1,2,\cdots,K,
  \end{array}
\end{equation}

where $\lambda_i(t) \in \mathbb{R}^{n_i}$ is the co-state vector, $\lambda=(\lambda_1^T,\lambda_2^T,\cdots \lambda_K^T)^T,$ and\\
$\Psi_i(x(t),\lambda(t))=\displaystyle \sum_{j=1}^{K}\frac{\partial f_j(x(t))}{\partial x_i(t)} \lambda_j(t).$ 
Also the optimal control law of the $ith$ subsystem is given by
\begin{equation}
  \label{Eq:Optimal control law}
  u_i^*(t)=-R_i^{-1}B_i^T\lambda_i(t), \quad t>t_0, \quad i=1,2,\cdots K.
\end{equation}

Unfortunately, problem (\ref{Eq:Optimality conditions 2}) is a nonlinear largescale TPBVP which is decomposed into N interconnected subproblems. In general, it is extremely
difficult to solve this problem analytically or even numerically, except in a few simple cases. In order to overcome this difficulty, we will presented the LaHOC method  in the next section.

\section{Laguerre polynomials and spectral homotopy analysis method }
In this section, we give a brief description of the basic idea of the Laguerre homotopy method for solving nonlinear boundary value problems. At first, we take into account the following properties of the modified Laguerre polynomials.

\subsection{Properties of the modified Laguerre polynomials}
Let $\omega_{\beta} (t)= e^{-\beta t}, \beta > 0,$ and define the weighted space ${L}_{\omega_{\beta}}^{2}(0,\infty)$  as usual, with the following inner product and norm, \cite{Ben}:
\begin{equation}
\label{C1}
(u, v)_{\omega_\beta}=\int_0^\infty u(t) v(t) \omega_{\beta} (t) dt, ~~~~ ||v||_{\omega_{\beta}}=(v,v)_{\omega_{\beta}}
\end{equation}
The modified Laguerre polynomial of degree l is defined by :
\begin{equation}
\label{C2}
\mathcal{L}_{l}^{\beta}(t)=\frac{1}{l!}e^{\beta t} \frac{d^l}{dt^l}(t^l e^{-\beta t}), ~~ l \geq 0.
\end{equation}
They satisfy the recurrence relation
\begin{equation}
\label{C3}
\frac{d}{dt}\mathcal{L}_{l}^{\beta}(t)= \frac{d}{dt} \mathcal{L}_{l-1}^{\beta}(t)-\beta \mathcal{L}_{l-1}^{\beta}(t), ~~ l \geq 1.
\end{equation}
The set of Laguerre polynomials is a complete ${L}_{\omega_{\beta}}^{2}(0,\infty)$-orthogonal system, namely,
\begin{equation}
\label{C4}
(\mathcal{L}_{l}^{\beta},\mathcal{L}_{m}^{\beta})_{\omega_{\beta}}=\frac{1}{\beta}\delta_{l,m},
\end{equation}
where $\delta_{l,m}$ is the Kronecker symbol. Thus, for any $v \in {L}_{\omega_{\beta}}^{2}(0,\infty),$
\begin{equation}
\label{C5}
 v(t)=\sum_{j=0}^{\infty}\hat{v}_l\mathcal{L}_{l}^{\beta}(t),
\end{equation}
where the coefficients $\hat{v}_l$ are given by
\begin{equation}
\label{C6}
\hat{v}_l=\beta(v,\mathcal{L}_{l}^{\beta})_{\omega_{\beta}}.
\end{equation}
Now, let $N$  be any positive integer, and $\mathcal{P}_{N}(0,\infty)$ the set of all algebraic polynomials of degree at most $N$. We denote by $t_{\beta,j}^N,\quad 0\leq j \leq N$ the nodes of modified Laguerre-Radau interpolation. Indeed, $t_{\beta,0}^N=0$ and $t_{\beta,j}^N,\quad 1\leq j \leq N$ are the distinct zeros
of  $\frac{d}{dt} \mathcal{L}_{N+1}^{\beta}(t)$ By using (\ref{C3}), the corresponding Christoffel numbers are as follows:
\begin{equation}
\label{C7}
\omega_{\beta,0}^N=\frac{1}{\beta(N+1)},~~~ \omega_{\beta,j}^N=\frac{1}{\beta(N+1)\mathcal{L}_{N}^{\beta}(t_{\beta,j}^N)~\mathcal{L}_{N+1}^{\beta}(t_{\beta,j}^N)},
\end{equation}
For any $\Phi \in  \mathcal{P}_{2N}(0,\infty)$,
\begin{equation}
\label{C8}
\sum_{j=0}^{N}\Phi(t_{\beta,j}^N) \omega_{\beta,j}^N=\int_0^\infty \Phi(t) \omega_{\beta} (t) dt.
\end{equation}

Next, we define the following discrete inner product and norm,
\begin{equation}
\label{C9}
(u, v)_{\omega_\beta,N}=\sum_{j=0}^{N}u(t_{\beta,j}^N)v(t_{\beta,j}^N)\omega_{\beta,j}^N, \qquad ||v||_{\omega_\beta,N}=(v,v)_{\omega_\beta,N}^{\frac{1}{2}}.
\end{equation}

For any $\Phi,\psi\in \mathcal{P}_{N}(0,\infty)$,
\begin{equation}
\label{C10a}
 (\Phi,\psi)_{\omega_\beta}=(\Phi,\psi)_{\omega_\beta,N}, \qquad ||v||_{\omega_\beta}=||v||_{\omega_\beta,N}.
\end{equation}


  \subsection{Spectral homotopy analysis method}
 In this section, we give a description of the SHAM with the Laguerre polynomials basis. This will be followed by a description of the new version of the SHAM  algorithm \cite{Motsa1}. To this end, consider a general $n$ dimensional initial value problem described as

\begin{align}
\mathbf{\dot{z}}(t) &= \mathbf{f}(t,\mathbf{z}(t)),\;\;\;\;\;\; \mathbf{z}(t_0) = \mathbf{z}^0, \label{seq1}\\
& \mathbf{z}:\mathbb{R} \rightarrow \mathbb{R} ^n,\;\;\;\; \mathbf{f}:\mathbb{R} \times \mathbb{R}^n \rightarrow \mathbb{R}^n
\end{align}
 We make the usual assumption that $\textbf{f}$ is sufficiently smooth for linearization techniques to be valid.  If $\mathbf{z} = (z_1,z_2,\ldots,z_n)$ we can apply the SHAM by rewriting equation (\ref{seq1}) as

\begin{equation}\label{seq2}
\dot{z}_r + \sum_{k=1}^{n} \sigma_{r,k} z_k + g_r(z_1,z_2,\ldots,z_n) = 0,
\end{equation}
subject to the initial conditions
\begin{equation}\label{seq3}
z_r(0) = z_r^0.
\end{equation}
where $z_r^0$ are the given initial conditions, $\sigma_{r,k}$ are known constant parameters and $g_r$ is the nonlinear component of the $r$th equation.

The SHAM approach imports the conventional ideas of the standard homotopy analysis method  by defining the following zeroth-order deformation equations
\begin{equation}\label{seq4}
(1-q)\mathcal{L}_r\left[\tilde{z}_r(t;q) - z_{r,0}(t) \right]
= q \hbar_r  \mathcal{N}_r[\tilde{\bf z}(t;q)],
\end{equation}
where $ q \in [0,1]$ is an embedding parameter, $\tilde{z}_r(t;q)$ are unknown functions, $\hbar_r$ is a
 convergence controlling parameter. The operators $\mathcal{L}_r$ and $\mathcal{N}_r$ are defined as
\begin{align}
 \mathcal{L}_r[\tilde{z}_r(t;q)] & =\frac{\partial \tilde{z}_r}{\partial t} + \sum_{k=1}^{n} \sigma_{r,k} \tilde{z}_k ,\label{seq5}\\
 \mathcal{N}_r[\tilde{\bf z}(t;q)] & = \mathcal{L}_r[\tilde{z}_r(t;q)] + g_r[\tilde{z}_1(t;q),\tilde{z}_2(t;q),\ldots, \tilde{z}_n(t;q)]\label{seq6}.
\end{align}

Using the ideas of the standard HAM approach \cite{Liao4}, we differentiate the zeroth-order equations (\ref{seq4}) $m$ times   with respect to $q$ and then set $q = 0$ and finally divide the resulting equations by $m!$ to obtain the following equations, which are referred to as the $m$th order  (or higher order) deformation equations,
\begin{equation}\label{seq7}
\mathcal{L}_r[z_{r,m}(t) - \chi_m z_{r,m-1}(t)]  = \hbar_r R_{r,m-1},
\quad m\ge 1,
\end{equation}
subject to
\begin{equation}\label{seq8}
z_{r,m}(0) = 0,
\end{equation}
where
\begin{equation}\label{seq9}
R_{r,m-1} = \left. \frac{1}{(m-1)!} \frac{\partial^{m-1} \mathcal{N}_r[\tilde{\bf z}(t;q)]}{\partial q^{m-1}} \right|_{q=0},
\end{equation}
and
\begin{equation}\label{seq10}
\chi_m=\left\{
\begin{array}{ll}
    0, & \hbox{$m\leqslant1$,} \\
    1, & \hbox{$m>1$.}\\
\end{array}%
\right.
\end{equation}

After obtaining solutions for equations (\ref{seq7}), the approximate solution for each $z_r(t)$ is determined as the series solution
\begin{equation}
\label{eq8}
 z_r(t)  = z_{r,0}(t) + z_{r,1}(t) + z_{r,2}(t) + \ldots
\end{equation}

A HAM solution is said to be of order $M$ if the above series is truncated at $m = M$, that is, if
\begin{equation}
\label{eq9}
 z_r(t)  = \sum_{m = 0}^{M} z_{r,m}(t).
\end{equation}
A suitable initial guess to start off the SHAM algorithm is obtained by solving the linear part of (\ref{seq2}) subject to the given initial conditions, that is, we solve

\begin{equation}
\label{eq10}
\mathcal{L}_r[z_{r,0}(t)]  = \phi_r(t),\;\;\;\;\; z_{r,0}(0) = z^0_{r}.
\end{equation}

If equation (\ref{eq10}) cannot be solved exactly, the spectral collocation method is used as a means of solution. The solution $ z_{r,0}(t)$ of equation (\ref{eq10}) is then fed to (\ref{seq7}) which is iteratively solved for $z_{r,m}(t)$ (for $m = 1,2,3 \ldots, M$).

In this paper, we use the Laguerre pseudo-spectral method to solve equations (\ref{seq7}-\ref{seq9}). The pseudo-spectral derivative $D_N(z)$ of a continuous function $z$ is defined by:
\begin{equation}
\label{C10a}
D_N(z)=D[I_N(z)],
\end{equation}
that is, $D_N(z)$ is the derivative of the interpolating polynomial of $z.$ Moreover, $D_N$ can be expressed in terms of a matrix, the pseudo-spectral derivation matrix $D_{\beta}$:
\begin{eqnarray*}
D_{\beta}=[(d_{\beta})_{ij}]_{i,j=0,1,\cdots,N}.
\end{eqnarray*}
Indeed, given the nodes $\{x_j^{(\beta)}\}_{j=0}^{N}$, an approximation $z\in \mathcal{P}^{(\beta)}_N$ of an unknown function and $\{(h_{\beta})_j\}$, 
the Lagrange interpolation polynomials associated to the points $x_j$, differentiating $m$ times the expression
\begin{eqnarray*}
	z_{\beta}(x)=\sum_{j=0}^{N} z_{\beta}(x_j)(h_{\beta})_j(x),
\end{eqnarray*}
yields:
\begin{eqnarray*}
z_{\beta}^{(m)}(x_k)=\sum_{j=0}^{N}(h_{\beta})_j^{(m)}(x_k)z_{\beta}(x_j),~~ 0\leq k \leq N.
\end{eqnarray*}
If we define:
\begin{eqnarray*}
&&z_{\beta}^{(m)} =\left(z_{\beta}^{(m)}(x_0), z_{\beta}^{(m)}(x_1),\cdots, z_{\beta}^{(m)}(x_N)\right)^T,~~ z_{\beta} = z_{\beta}^{(0)},\\
&& D_{\beta}^{(m)} =\left[(d_{\beta})_{ij}^{(m)}=(h_{\beta})_{j}^{(m)}(x_i)\right]_{0\leq i,j \leq N},\\
&& (d_{\beta})_{ij}^{(m)}=(h_{\beta}
)_{j}^{(m)}(x_i),
\end{eqnarray*}
then:
$$D_{\beta} = D_{\beta}^{(1)}, (d_{\beta})_{ij} = (d_{\beta})_{ij}^{(1)}.$$
We now state two important results. The first ensures that it is sufficient to compute the first
order differentiation matrix, the second gives the general expression of its entries.
\begin{lem}\label{Theo1}\cite{JieShen}
\begin{equation}
D_{\beta}^{(m)}=D_{\beta}.D_{\beta}\cdots D_{\beta}=D_{\beta}^m,~~m\geq 1.
\end{equation}
\end{lem}

Let $\{x_j^{(\beta)}\}_{j=0}^{N}$ be the Gauss-Laguerre (GL) or Gauss-Laguerre-Radau (GLR) nodes and $z\in  \mathcal{P}^{(\beta)}_N.$ Let $\{(h_{\beta})_j(x)\}_{j=0}^N$ be the Lagrange interpolation polynomials relative to $\{x_j^{(\beta)}\}_{j=0}^{N}$. From  Lemma \ref{Theo1}, we have: 
\begin{eqnarray*}
z_{\beta}^{(m)} = D_{\beta}^m z_{\beta},~ m\geq 1.
\end{eqnarray*}
Next we have:
\begin{lem} \cite{JieShen}
The entries of the differentiation matrix $D_{\beta}$ associated to the GL and GLR points $\{x_j^{(\beta)}\}_{j=0}^N$ have the following form:
\begin{itemize}
\item GL points: $\{x_j^{(\beta)}\}_{j=0}^N$ are the zeros of $\mathscr{L}^{(\beta)}_{N+1}(x)$, 
\begin{equation}
\label{eq11}
d_{ij} =
\left\{\begin{array}{ll}
\frac{\mathscr{L}_N^{(\beta)}\left(x_i^{(\beta)}\right)}{\left(x_i^{(\beta)}-x_j^{(\beta)}\right)\mathscr{L}_N^{(\beta)}\left(x_j^{(\beta)}\right)} & \textrm{if}\;\; i\neq j,\\
\\
 \frac{\beta x_i^{(\beta)}-N-2}{2 x_i^{(\beta)}}&\textrm{if} \;\;  i =j,
\end{array}\right.
\end{equation}

\item GLR points: $x_0=0$, $\{x_j^{(\beta)}\}_{j=1}^N$ are the zeros of $\frac{\partial}{\partial x} \mathscr{L}_{N+1}^{(\beta)}(x)$,
\begin{equation}
\label{eq12}
d_{ij} =
\left\{\begin{array}{ll}
\frac{\mathscr{L}_{N+1}^{(\beta)}\left(x_i^{(\beta)}\right)}{\left(x_i^{(\beta)}-x_j^{(\beta)}\right)\mathscr{L}_{N+1}^{(\beta)}\left(x_j^{(\beta)}\right)} & \textrm{if}\;\; i\neq j,\\
\\
 \frac{\beta}{2}&\textrm{if} \;\;  i =j\neq 0,\\
 \\
\frac{-\beta N}{2}&\textrm{if} \;\;  i =j=0, 
 \end{array}\right.
\end{equation}
\end{itemize}
\end{lem}

Applying the the Laguerre spectral collocation method in equations (\ref{seq7}-\ref{seq9}) gives

\begin{equation}
\label{eq13}
\textbf{A}\left[ \textbf{W}_{m} -  \chi_m\textbf{W}_{m-1} \right] =   \hbar_r \textbf{R}_{m-1},  \;\;\;\;\;\;\; \textbf{W}_{m}(\tau_0)=0, \;\;\;\;\;\;\; \textbf{W}_{m}(\tau_N)=0,
\end{equation}
where   $\textbf{R}_{m-1}$ is an  $(N+1)n\times 1$ vector corresponding to $R_{r,m-1}$ when evaluated at the collocation points and $\textbf{W}_m = [\tilde{\textbf{z}}_{1,m}; \tilde{\textbf{z}}_{2,m} ;\ldots ; \tilde{\textbf{z}}_{n,m}]$.

The matrix $\textbf{A}$ is an $(N+1)n\times(N+1)n$ matrix that is derived from transforming the linear operator $\mathcal{L}_r$  using the derivative matrix $D_{\beta} $ (we omit subscipt $\beta$ for simplicity) and is defined as

\begin{equation}
\label{eq14}
\textbf{A}  = \left[\begin{array}{cccc}
A_{11}  & A_{12} & \cdots & A_{1n}\\
A_{21}  & A_{22} & \cdots & A_{2n}\\
 \vdots &        &  \ddots & \vdots \\
A_{n1}  & A_{n2}  & \cdots & A_{nn}
\end{array}\right],\;\;\; \textrm{with}\;\;\;\; \textbf{A}_{pq} = \left\{\begin{array}{cr}
\displaystyle  \textbf{D} + \sigma_{pq} \textbf{I} ,& p = q, \\
\displaystyle \sigma_{pq} \textbf{I} ,&  p \neq q,
\end{array} \right.
\end{equation}
where $\textbf{I}$ is an identity matrix of order $N+1$.

Thus, starting from the initial approximation, the recurrence formula (\ref{eq13}) can be used to obtain the solution $z_r(t)$.\\

\section{Convergence analysis of LaHOC}
To analysis the convergence of LaHOC, we first recall the  $m$th order  (or higher order) deformation equation,
\begin{equation}
\label{eq1H}
\mathcal{L}[z_{m}(t) - \chi_m z_{m-1}(t)]  = \hbar H(t)R_{m-1},
\end{equation}
 subject to the initial condition
\begin{equation}
\label{eq2H}
  z_{m,1:n}(t_0) = 0,
\end{equation}
 where  $H(t)\neq 0$  is an auxiliary function, 
\begin{eqnarray}
\label{eq3H}
R_{m-1}=\mathcal{L}[z_{m-1}]+\mathcal{N}_{m-1}[z_0,z_1,\cdots, z_{m-1}]-(1-\chi_{m})\phi(t).
\end{eqnarray}
where $z_{r,m}$, $\mathcal{L}_r$ and $\mathcal{N}_r$ in (\ref{seq7}) are the $r$th components of $z_{m-1}$ and operators $\mathcal{L}$
and $\mathcal{N}$, respectively.
Let us define the nonlinear operator $\mathcal{N}$ and the sequence $\{Z_{m}\}_{m=0}^{\infty}$ as,
\begin{eqnarray}
\mathcal{N}[\mathbf{z}(t)]=\sum_{k=0}^{\infty}N_k(z_0,z_1,\cdots, z_k),
\end{eqnarray}
\begin{equation}
\label{eq4H}
\left\{ \begin{array}{ll}
Z_0=z_0,\\
Z_1=z_0+z_1,\\
\vdots\\
Z_m=z_0+z_1+z_2+\cdots+ z_m.
\end{array}\right.
\end{equation}
Therefore,  we have
\begin{eqnarray}
\label{eq7a}
\mathcal{L}[z_{m}(t)] = \hbar H(t)\{\sum_{k=0}^{m-1}\mathcal{L}[z_{k}]+\sum_{k=0}^{m-1}\mathcal{N}_{k}-\phi(t)\},
\end{eqnarray}
from (\ref{eq4H}) we have
\begin{eqnarray}
\label{eq8a}
\mathcal{L}[Z_{m}(t)-Z_{m -1}(t)] = \hbar H(t)\{\mathcal{L}[Z_{m-1}]+\mathcal{N}[Z_{m-1}]-\phi(t)\},
\end{eqnarray}
subject to the initial condition
\begin{eqnarray}
\label{eq9a}
Z_{m,1:n}(t_0) = 0.
\end{eqnarray}
Consequently, the collocation method is based on a solution $Z^N(t)\in \mathcal{P}_{N+1}(0, \infty)$, for (\ref{eq8a}) such that
\begin{eqnarray}
\label{eq13a}
\mathcal{L}[Z^N_{m}(t_{\beta,k}^N)-Z^N_{m-1}(t_{\beta,k}^N)] = \hbar H^N(t_{\beta,k}^N)\{\mathcal{L}[Z^N_{m-1}(t_{\beta,k}^N)]+\mathcal{N}[Z^N_{m-1}(t_{\beta,k}^N)]-\phi^N(t_{\beta,k}^N)\},
\end{eqnarray}
subject to the initial condition
\begin{eqnarray}
\label{eq14a}
Z^N_{m,1:n}(t_0) = 0.
\end{eqnarray}
From (\ref{eq13a}) we have
\begin{eqnarray}
\label{eq16a}
&&\mathcal{L}[Z^N_{m}(t_{\beta,k}^N)] = (1+\hbar H^N(t_{\beta,k}^N))\mathcal{L}[Z^N_{m-1}(t_{\beta,k}^N)]+\hbar H^N(t_{\beta,k}^N)\{\mathcal{N}[Z^N_{m-1}(t_{\beta,k}^N)]-\phi^N(t_{\beta,k}^N)\},\nonumber\\
&&~~~~~~~~~~~~~~~~\quad 0\leq k \leq N,\quad  m\geq1,\\
&&Z^N_{m,1:n}(t_0) = 0, m\geq 0.\nonumber
\end{eqnarray}
Now, we choose $L[Z(t)]=\frac{\mathrm{d}}{\mathrm{d}t}Z+\alpha(t)Z$, $N[Z(t)]=-\alpha(t)Z-f(t,Z)$ and $\phi(t)\equiv0$ where $\alpha(t)$ is an arbitrary analytic function.\\
Let $\tilde{Z}^{N}_{m}(t) = Z^{N}_{m}(t)-Z^{N}_{m-1}(t) $, then we have from (\ref{eq16a}) that
\begin{eqnarray}
\label{eq17a}
&&\mathcal{L}[\tilde{Z}^{N}_{m}(t_{\beta,k}^N)]=(1+\hbar H(t_{\beta,k}^N))\mathcal{L}[Z^{N}_{m-1}(t_{\beta,k}^N)-Z^{N}_{m-2}(t_{\beta,k}^N)]+\hbar H(t_{T,k}^N)\nonumber\\
&&~~~~~~~~~~~~~~~~~~~~\{\mathcal{N}[Z^{N}_{m-1}(t_{\beta,k}^N)]-\mathcal{N}[Z^{N}_{m-2}(t_{\beta,k}^N)]\},\nonumber\\
&&~~~~~~~~~~~~~~~~~~~~~~~~~~~~~~\quad 0\leq k \leq N,  m\geq1,
\end{eqnarray}
or according to the definitions of $L[Z(t)]$ and $N[Z(t)]$,
\begin{eqnarray}
\label{eq17b}
&&\frac{d}{dt}[\tilde{Z}^{N}_{m}(t_{\beta,k}^N)]+\alpha(t_{\beta,k}^N)\tilde{Z}^{N}_{m}=(1+\hbar H(t_{\beta,k}^N))\frac{d}{dt}[\tilde{Z}^{N}_{m-1}(t_{\beta,k}^N)]+\alpha(t_{\beta,k}^N)\tilde{Z}^{N}_{m-1} \nonumber\\
&&~~~~~~~~~~~~~~~~~~~~~~~~~~~~~~-\hbar H(t_{\beta,k}^N)
\{{f}(t_{\beta,k}^N,Z^{N}_{m-1}(t_{\beta,k}^N)) -{f}(t_{\beta,k}^N,Z^{N}_{m-2}(t_{\beta,k}^N))\},\nonumber\\
&&~~~~~~~~~~~~~~~~~~~~~~~~~~~~~~\quad 0\leq k \leq N,  m\geq1,
\end{eqnarray}

\begin{theorem} Assume that for any $k = 0, 1, . . . ,N, \mathscr{Z}_k=\{Z^{N}_{m}(t_{\beta,k}^N)\}_0^\infty$ is the LaHOC sequence produced by (\ref{eq16a}). Furthermore, assume
$\alpha_0=\min_{\,t\in[0,\infty)}\alpha(t)$,
$\alpha_1=\max_{\,t\in[0,\infty)}|\alpha(t)|$ and $H=\max_{\,t\in[0,\infty)}|H(t)|$ and
\begin{eqnarray}
\label{eqlip}
||f(., Z^{N}_{m})-f(., Z^{N}_{m-1})||_{\omega_\beta,N}\leq L_f ||Z^{N}_{m}-Z^{N}_{m-1}||_{\omega_\beta,N}.
\end{eqnarray}
for some constant $L_f > 0$. Then for any initial n-vector $Z^{N}_{0}(t_{\beta,k}^N),$ $\mathscr{Z}_k$ converges to some  $\hat{Z}(t_{\beta,k}^N)$ which is the exact solution of  (\ref{seq2}), at any GLR point, $t_{\beta,k}^N$, if
\begin{eqnarray}
\label{Cr1}
\gamma=\frac{N|1+ \hbar H| + \alpha_1 + |\hbar| H L_f }{\beta/2+\alpha_0}<1.
\end{eqnarray}
\end{theorem}

\noindent\textbf{Proof.}
1. Using (\ref{C1}) and integrating by parts yield that
\begin{eqnarray}
\label{eq18a}
\left(\tilde{Z}^{N}_{m},\frac{d}{dt}\tilde{Z}^{N}_{m} \right)_{\omega_{\beta},N}=\left(\tilde{Z}^{N}_{m},\frac{d}{dt}\tilde{Z}^{N}_{m} \right)_{\omega_{\beta}}=\frac{1}{2}\left[e^{-\beta t}(\tilde{Z}^{N}_{m})^{2}\mid _0^\infty+\int_{0}^{\infty}\beta e^{-\beta t}(\tilde{Z}^{N}_{m})^{2}dt\right],
\end{eqnarray}
then, we have
\begin{eqnarray}
\label{eq19a}
2\left(\tilde{Z}^{N}_{m},\frac{d}{dt}\tilde{Z}^{N}_{m} \right)_{\omega_{\beta},N}=\beta\|\tilde{Z}^{N}_{m}\|_{\omega_\beta}^{2},~~~~~\|\tilde{Z}^{N}_{m}\|_{\omega_\beta,N}=\|\tilde{Z}^{N}_{m}\|_{\omega_\beta},
\end{eqnarray}
by (\ref{eq19a}) and  from the Cauchy inequality we obtain that
\begin{eqnarray}
\label{eq20a}
\beta\|\tilde{Z}^{N}_{m}\|_{\omega_\beta}^{2}\leq 2 \|\tilde{Z}^{N}_{m}\|_{\omega_\beta,N}\| \frac{d}{dt}(\tilde{Z}^{N}_{m})\|_{\omega_\beta,N},
\end{eqnarray}
from where
\begin{eqnarray}
\label{eq21a}
\|\tilde{Z}^{N}_{m}\|_{\omega_\beta}\leq \frac{2}{\beta}\|\frac{d}{dt}(\tilde{Z}^{N}_{m})\|_{\omega_\beta},
\end{eqnarray}

2. Taking discrete weighted inner product of (\ref{eq17b}) with $\tilde{Z}^{N}_{m}(t_{\beta,k}^N)$, we have
\begin{eqnarray}
\label{eq17c}
&&\left(\frac{d}{dt}\tilde{Z}^{N}_{m}+\alpha(t)\tilde{Z}^{N}_{m}, \tilde{Z}^{N}_{m} \right)_{\omega_\beta, N}
=\left((1+\hbar H) \frac{d}{dt}\tilde{Z}^{N}_{m-1}+\alpha(t)\tilde{Z}^{N}_{m-1}
,\tilde{Z}^{N}_{m}\right)_{\omega_\beta, N} \nonumber\\
&&~~~~~~~~~~~~~~~~~~~~~~~~~~~~~~-\hbar 
\left(H(t)[f(t_{\beta,k}^N,Z^{N}_{m-1} -f(t_{\beta,k}^N,Z^{N}_{m-2})],
\tilde{Z}^{N}_{m}\right)_{\omega_\beta, N}\nonumber\\
&&~~~~~~~~~~~~~~~~~~~~~~~~~~~~~~\quad 0\leq k \leq N,  m\geq1,
\end{eqnarray}

Therefore, a combination with Cauchy inequality and \eqref{eq19a} leads to
\begin{eqnarray}
\label{eq22a}
(\frac\beta 2+\alpha_0)\|\tilde{Z}^{N}_{m}\|_{\omega_\beta} & \leq &
|1+ \hbar H| \|\frac{d}{dt}\tilde{Z}^{N}_{m-1}\|_{\omega_\beta}
+ \alpha_1 ||\tilde{Z}^{N}_{m-1}||_{\omega_\beta}\nonumber\\
&& + |\hbar| H \| f(t_{\beta,k}^N,Z^{N}_{m-1} -f(t_{\beta,k}^N,Z^{N}_{m-2}) \|_{\omega_\beta,N} 
\end{eqnarray}
Then by using inverse inequality of Laguerre polynomial and \eqref{eqlip}, we get
\begin{equation}
\label{eq23a}
(\frac\beta 2+\alpha_0)\|\tilde{Z}^{N}_{m}\|_{\omega_\beta} \leq 
(N|1+ \hbar H| + \alpha_1 + |\hbar| H L_f ) \|\tilde{Z}^{N}_{m-1}\|_{\omega_\beta},
\end{equation}
%
%
which is
\begin{equation}
\label{25a}
    \left\|\tilde{Z}^{N}_{m}\right\|_{\omega_\beta}\leq  
    \frac{N|1+ \hbar H| + \alpha_1 + |\hbar| H L_f }{\beta/2+\alpha_0}
    \left\|\tilde{Z}^{N}_{m-1}\right\|_{\omega_\beta} = \gamma \left\|\tilde{Z}^{N}_{m-1}\right\|_{\omega_\beta}.
\end{equation}
Hence, we have
\begin{equation}
\label{25b}
    \left\|\tilde{Z}^{N}_{m}\right\|_{\omega_\beta}\leq \gamma \left\|\tilde{Z}^{N}_{m-1}\right\|_{\omega_\beta}\leq \cdots \leq \gamma ^m\left\|\tilde{Z}^{N}_{0}\right\|_{\omega_\beta}.
\end{equation}
Then for any $m' \geq m \geq 1,$
\begin{equation}
\label{25b}
\left\|{Z}^{N}_{m'}-{Z}^{N}_{m}\right\|_{\omega_\beta}\leq \sum_{i=m+1}^{m'}\left\|\tilde{Z}^{N}_{i}\right\|_{\omega_\beta}\leq \sum_{i=m+1}^{m'}\gamma^i \left\|\tilde{Z}^{N}_{0}\right\|_{\omega_\beta} \leq \frac{\gamma^{m+1}}{1-\gamma} \left\|\tilde{Z}^{N}_{0}\right\|_{\omega_\beta}.
\end{equation}
Since $\gamma \in [0, 1), \left\|{Z}^{N}_{m'}-{Z}^{N}_{m}\right\|_{\omega_\beta} \rightarrow 0$ as $m, m' \rightarrow \infty$. Thus $\mathscr{Z}_k$ is a Cauchy sequence; and since $\mathbb{R}^n$ is a Banach space, $\mathscr{Z}_k$ has a limit $\hat{Z}(t_{\beta,k}^N)$. Taking limit $m \rightarrow \infty$ in {\color{red}\eqref{eq13a}}, yields
\begin{eqnarray*}
\label{25b}
&\mathcal{L}[\hat{Z}(t_{\beta,k}^N)-\hat{Z}(t_{\beta,k}^N)] =0= \hbar H(t_{\beta,k}^N)\{\mathcal{L}[\hat{Z}(t_{\beta,k}^N)]+\mathcal{N}[\hat{Z}(t_{\beta,k}^N)]-\phi^N(t_{\beta,k}^N)\},\\
& \hat{Z}(0) = z^0.
\end{eqnarray*}
Thus, $\hat{Z}(t_{\beta,k}^N)$ is the exact solution of (\ref{seq2})  at any GLR point $t_{\beta,k}^N$. Also, by noticing the definition of $\hat{Z}^N(t),$ it is easy to verify $\hat{Z}^N(t_{\beta,k}^N)=\hat{Z}(t_{\beta,k}^N)$ and the proof is completed.


 \section{Numerical experiments}
 To demonstrate the applicability of the LaHOC algorithm as an
 appropriate tool for solving  infinite horizon optimal control for nonlinear large-scale dynamical systems, we apply the proposed algorithm to several test problems.\\

 \noindent\textbf{ Test problem 3.1.} Consider the two-order nonlinear composite system described by \cite{Tang1}:\\
  \begin{eqnarray}
    \label{system1}
    &&\dot x_1(t)=x_1(t)+u_1(t)-x_1^3(t)+x_2^2(t),\\
    &&\dot x_2(t)=-x_2(t)+u_2(t)+x_1(t)x_2(t)+x_2^3(t),\\
    &&x_1(0)=0, \quad x_2(0)=0.8.
  \end{eqnarray}
  The quadratic cost functional to be minimized is given by:
  \begin{eqnarray}
    \label{Performance index m}
    J=\frac{1}{2}\sum_{i=1}^{2}\int_{0}^{\infty}(x_i^2(t)+u_i^2(t))dt,
  \end{eqnarray}
  In this example, we have  $A_1=B_1=B_2=1$,  $A_2=-1$,  $Q_1=Q_2=R_1=R_2=1$,  $f_1(x)=-x_1^3(t)+x_2^2(t)$, $f_2(x)=x_1(t)x_2(t)+x_2^3(t)$.\\
  Then, according to the optimal control theory
  (\ref{Eq:Optimality conditions 2}), the optimality
  conditions can be written as:
  \begin{eqnarray}
    \label{system1}
    &&\dot x_1(t)=x_1(t)-\lambda_1(t)-x_1^3(t)+x_2^2(t),\\
    &&\dot x_2(t)=-x_2(t)-\lambda_2(t)+x_1(t)x_2(t)+x_2^3(t),\\
    &&\dot \lambda_1(t)=-x_1(t)-\lambda_1(t)+3x_1^2(t)\lambda_1(t)-x_2(t)\lambda_2(t),\\
    &&\dot \lambda_2(t)=-x_2(t)+\lambda_2(t)-2x_2(t)\lambda_1(t)-x_1(t)\lambda_2(t)-3x_2^2(t)\lambda_2(t),\\
    &&x_1(0)=0, \quad x_2(0)=0.8,\quad \lambda_1(\infty)=0, \quad \lambda_2(\infty)=0.
  \end{eqnarray}
 Also the optimal control laws are $u_1(t)=-\lambda_1,~   u_2(t)=-\lambda_2.$ 
  
  In this example, the parameters used in the LaHOC algorithms are

  \begin{align}
    &\mathcal{L}_r = \left[\begin{array}{cccc}
        \frac{d}{dt}-1 & 0 & 1 & 0 \\
        0 & \frac{d}{dt}+1 & 0 & 1 \\
        1   &  0 & \frac{d}{dt}+1 & 0 \\
        0   &  0 & 0 & \frac{d}{dt}-1 \\
      \end{array}\right],\;
    \textbf{A} = \left[\begin{array}{cccc}
        \textbf{D}-I& O & I  & O \\
        O &  \textbf{D} +I & O & I\\
        I & O &  \textbf{D}+I & O \\
        O & O &  O & \textbf{D}-I \\
      \end{array}\right],
  \end{align}
  \begin{align}
    &\mathcal{F}_r = \left[\begin{array}{c} x_1^3-x_2^2 \\
       - x_1 x_2-x_2^3\\-3x_1^2 \lambda_1 +x_2 \lambda_2
        \\2x_2\lambda_1+x_1\lambda_2+3x_2^2\lambda_2
      \end{array}\right],\; \phi = \left[\begin{array}{c}
        0 \\ 0 \\0\\ 0\\0\\
      \end{array}\right],\\ &R_{r,m-1} = \mathcal{L}_r[x_{r,m-1}] + Q_{r,m-1},\\ 
    &Q_{r,m-1} = \left[\begin{array}{l}
        -\displaystyle\sum_{j=0}^{m-1}\textbf{Z}_{1,m-1-j}(t)\displaystyle\sum_{k=0}^{j}\textbf{Z}_{1,j}(t)\textbf{Z}_{1,j-k}(t)+\displaystyle \sum_{j=0}^{m-1}\textbf{Z}_{2,j} \textbf{Z}_{2,m-1-j}\\
        \displaystyle \sum_{j=0}^{m-1}\textbf{Z}_{1,j}(t) \textbf{Z}_{2,m-1-j}(t)+\displaystyle\sum_{j=0}^{m-1}\textbf{Z}_{2,m-1-j}(t)\displaystyle\sum_{k=0}^{j}\textbf{Z}_{2,j}(t)\textbf{Z}_{2,j-k}(t)\\
        3\displaystyle\sum_{j=0}^{m-1}\textbf{Z}_{3,m-1-j}(t)\displaystyle\sum_{k=0}^{j}\textbf{Z}_{1,j}(t)\textbf{Z}_{1,j-k}(t)-\displaystyle \sum_{j=0}^{m-1}\textbf{Z}_{2,j}(t) \textbf{Z}_{4,m-1-j}(t)\\
        -2\displaystyle \sum_{j=0}^{m-1}\textbf{Z}_{2,j}(t) \textbf{Z}_{3,m-1-j}(t)-\displaystyle \sum_{j=0}^{m-1}\textbf{Z}_{1,j}(t) \textbf{Z}_{4,m-1-j}(t)-3\displaystyle\sum_{j=0}^{m-1}\textbf{Z}_{4,m-1-j}(t)\displaystyle\sum_{k=0}^{j}\textbf{Z}_{2,j}(t)\textbf{Z}_{2,j-k}(t)\\
      \end{array}\right]
  \end{align}
  With these definitions, the LaHOC algorithm gives
  \begin{equation}\label{soln1}
    \textbf{X}_{r,m} =  (\chi_m + \hbar_r )\textbf{X}_{r,m-1}  + \hbar_r \textbf{A}^{-1}\textbf{Q}_{r,m-1},
  \end{equation}
  Because the right hand side of equation (\ref{soln1}) is known, the solution can easily be obtained by using methods for solving linear system of equations.

Table 1 gives a comparison between the present LaHOC results for  $N = 100$ and $\hbar = -0.6$ and the numerically generated $\mathtt{BVP5C}$  \cite{KierS2008},  at selected values of time $t.$ It can be seen from the table that there is in good agreement between the two results.  
Moreover, our calculations show the better accuracy of LaHOC. In comparison with the $\mathtt{BVP5C}$, it is noteworthy that the LaHOC controls the error bounds while preserving the CPU time. The CPU time of LaHOC is $0.606532~s$, and $\mathtt{BVP5C}$ is $1.109817~s.$

Figure 1 and Figure 2 show the suboptimal states and control for $m = 19$  iterations
of LaHOC, compared to  MATLAB built-in function $\mathtt{BVP5C}$. The convergence of LaHOC ieteration is depicted in Figure 3. Also, Figure 4 presents that the minimum objective functional $|J_{j}-J_{N}|$ converges to $0$, where $j = 20,30,\ldots, 110$ and $N = 120.$

The results obtained with the present method are in good agreement with results of the successive approximation method used
by Tang and Sun \cite{Tang1}. 

       
\begin{table}[!h]
 \caption{Comparison between the LaHOC solution when $N = 100$ and $\hbar = -0.6$  and $\mathtt{BVP5C}$  solution.}
       \label{Label}\centering
      \begin{tabular*}{0.99\textwidth}{@{\extracolsep{\fill}}ccccccccc}\\
	\hline
     &  & $x_1(t)$  &  & $x_2(t)$ & & $\lambda_1(t)$ & & $\lambda_2(t)$ \\
     \cline{2-3}\cline{4-5}\cline{6-7}\cline{8-9}
       \\
  \hline   
  $t$    & $LaHOC$ & $BVP5C$ & $LaHOC$ & $BVP5C$ & $LaHOC$ & $BVP5C$ & $LaHOC$ & $BVP5C$\\
        \hline
0.113  &	 0.013872  &	 0.013872  &	 0.689067  &	 0.689067 &	 0.388387  &	 0.388387  &	 0.557556  &	 0.557556   \\
0.494  &	 0.031434  &	 0.031434  &	 0.412872  &	 0.412872 &	 0.195039  &	 0.195039  &	 0.236820  &	 0.236820   \\
1.152  &	 0.021573  &	 0.021573  &	 0.164529  &	 0.164529 &	 0.070317  &	 0.070317  &	 0.075704  &	 0.075704   \\
2.107  &	 0.006800  &	 0.006800  &	 0.042594  &	 0.042594 &	 0.017627  &	 0.017627  &	 0.018077  &	 0.018077   \\
3.389  &	 0.001168  &	 0.001168  &	 0.006943  &	 0.006943 &	 0.002852  &	 0.002852  &	 0.002887  &	 0.002887   \\
5.047  &	 0.000113  &	 0.000113  &	 0.000666  &	 0.000666 &	 0.000273  &	 0.000273  &	 0.000276  &	 0.000276   \\
 \hline
   \end{tabular*}
       \end{table}
  
  \noindent\begin{tabular*}{.5cm}{cc}
 \includegraphics[width=0.480\textwidth]{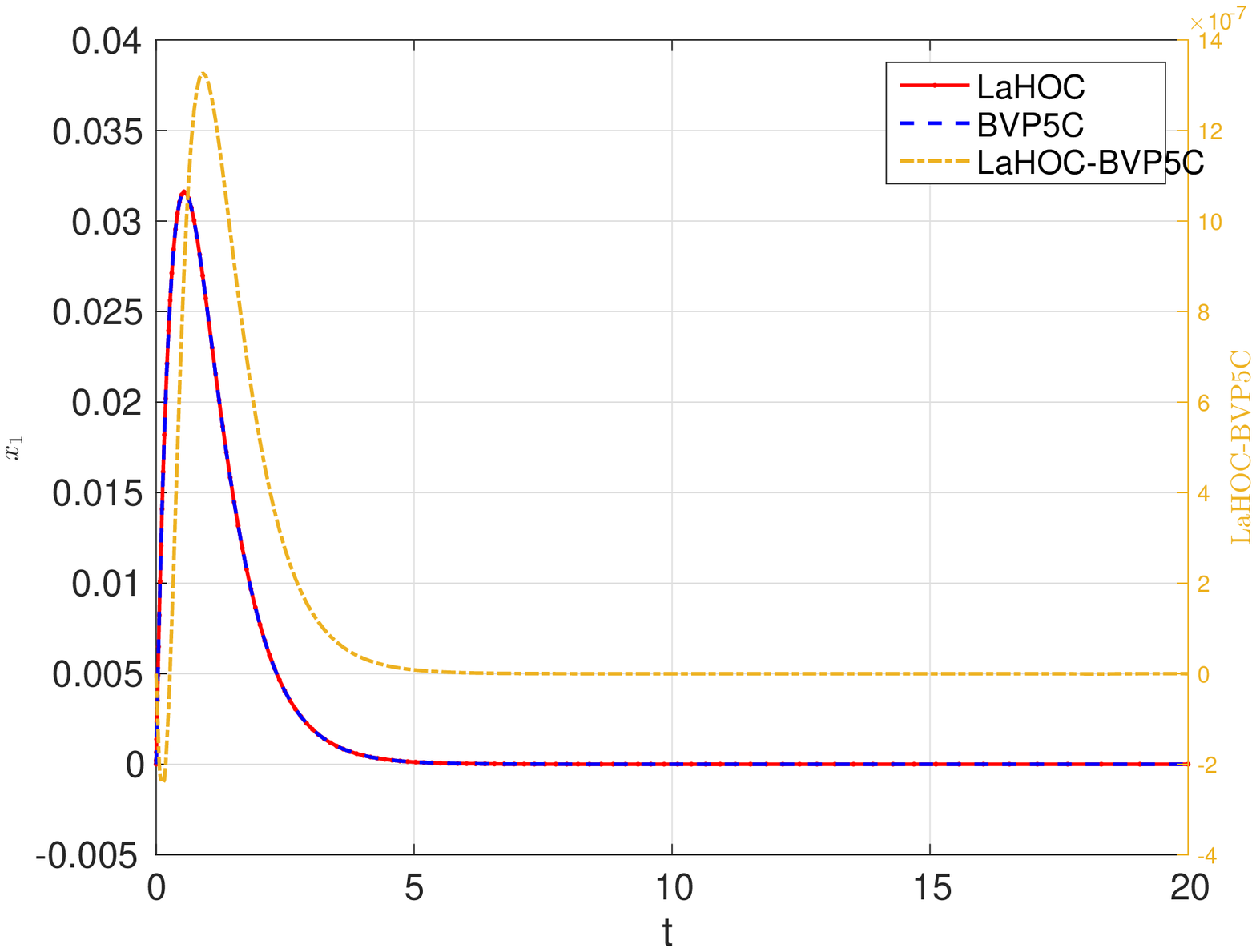} \includegraphics[width=0.480\textwidth]{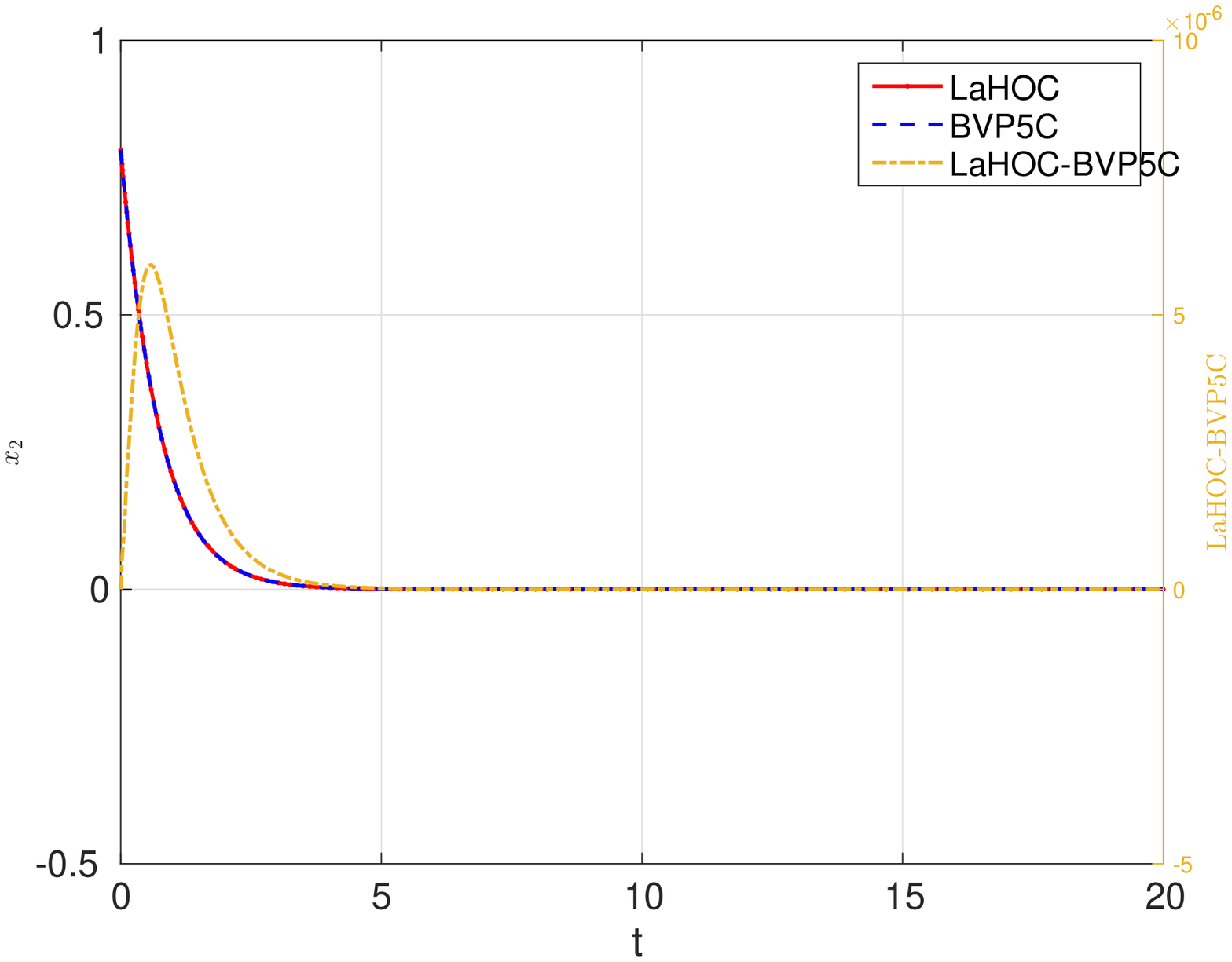}\\
 \centerline{Fig. 1.  The amplitudes of optimal state variables.}  \\
 \end{tabular*}\\
\begin{tabular*}{.5cm}{cc}
 \includegraphics[width=0.480\textwidth]{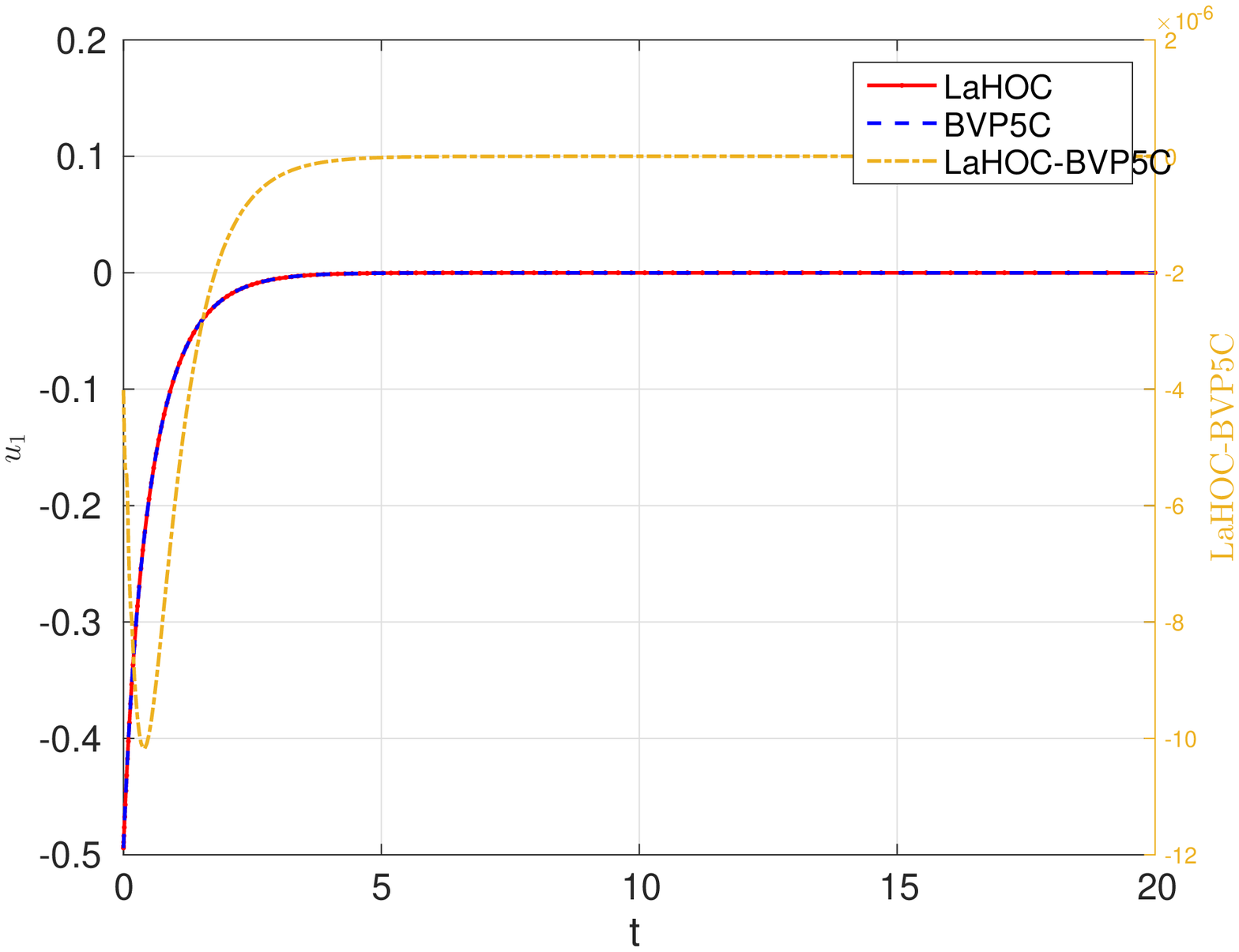} \includegraphics[width=0.480\textwidth]{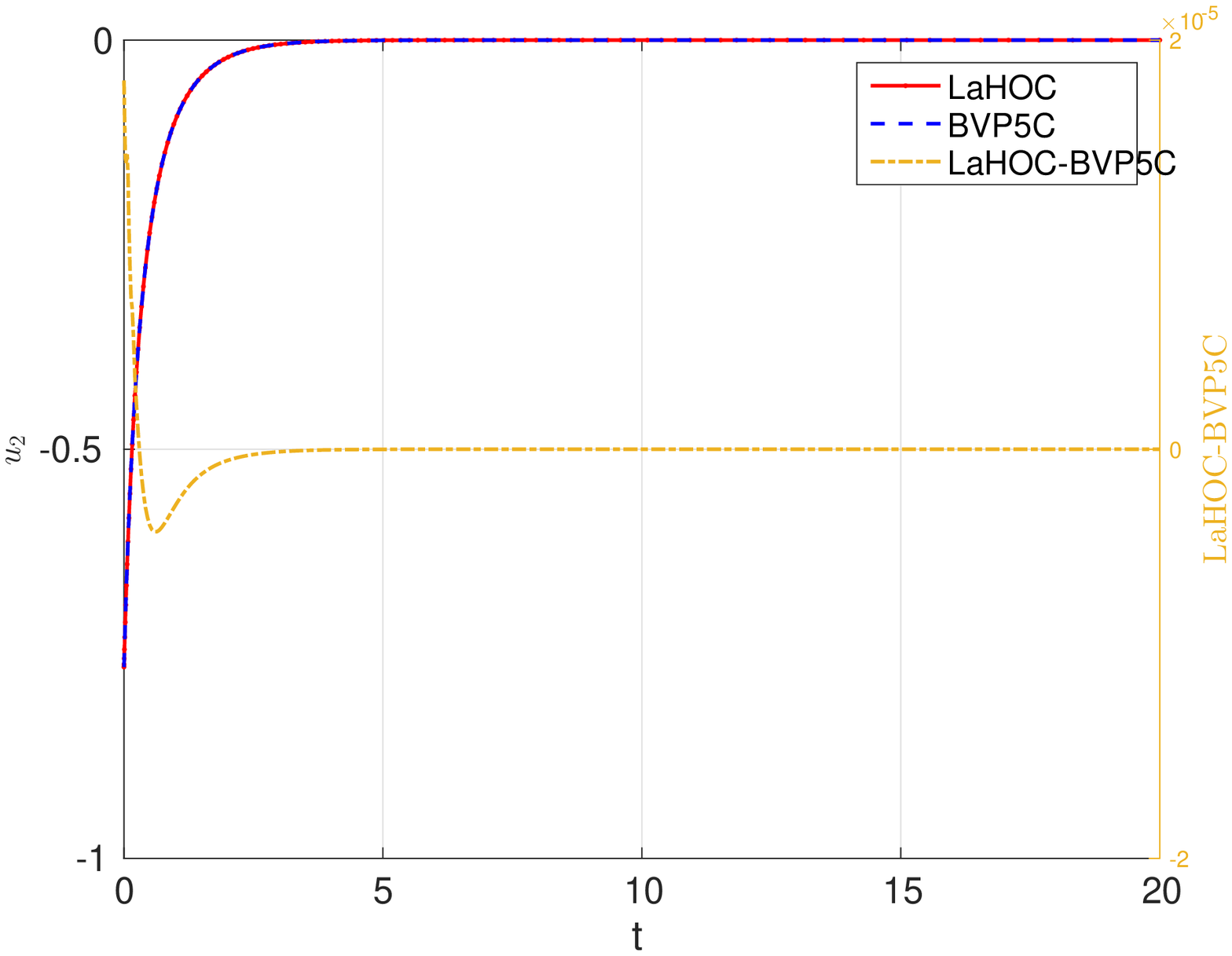}\\
 \centerline{Fig. 2.  The amplitudes of optimal control variables.}  \\
 \end{tabular*}\\
 \begin{tabular*}{.5cm}{cc}
  \includegraphics[width=0.50\textwidth]{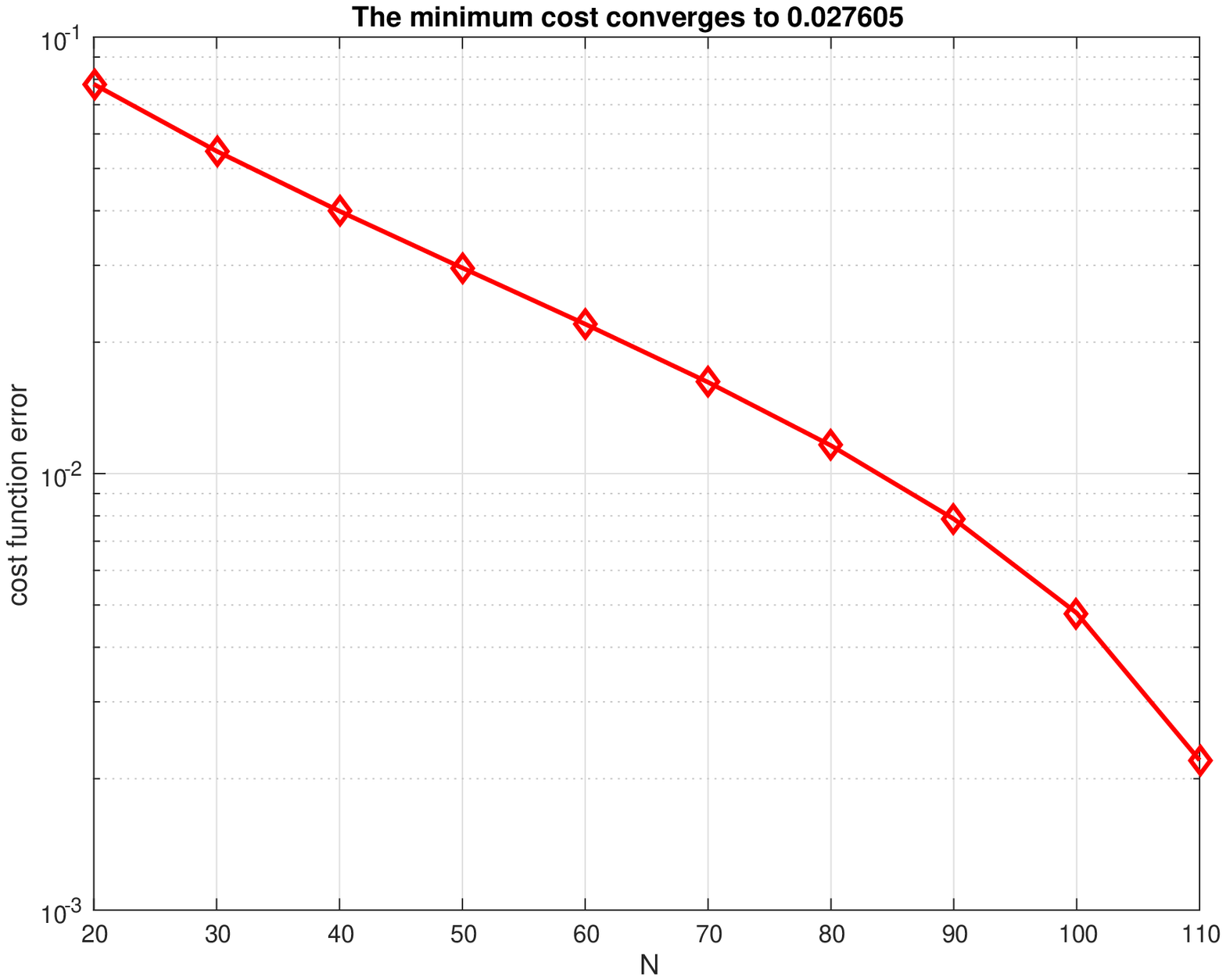}\\
 \centerline{Fig. 3.  The minimum cost onvergence.}  \\
 \end{tabular*}\\
\begin{tabular*}{.5cm}{cc}
 \includegraphics[width=0.450\textwidth]{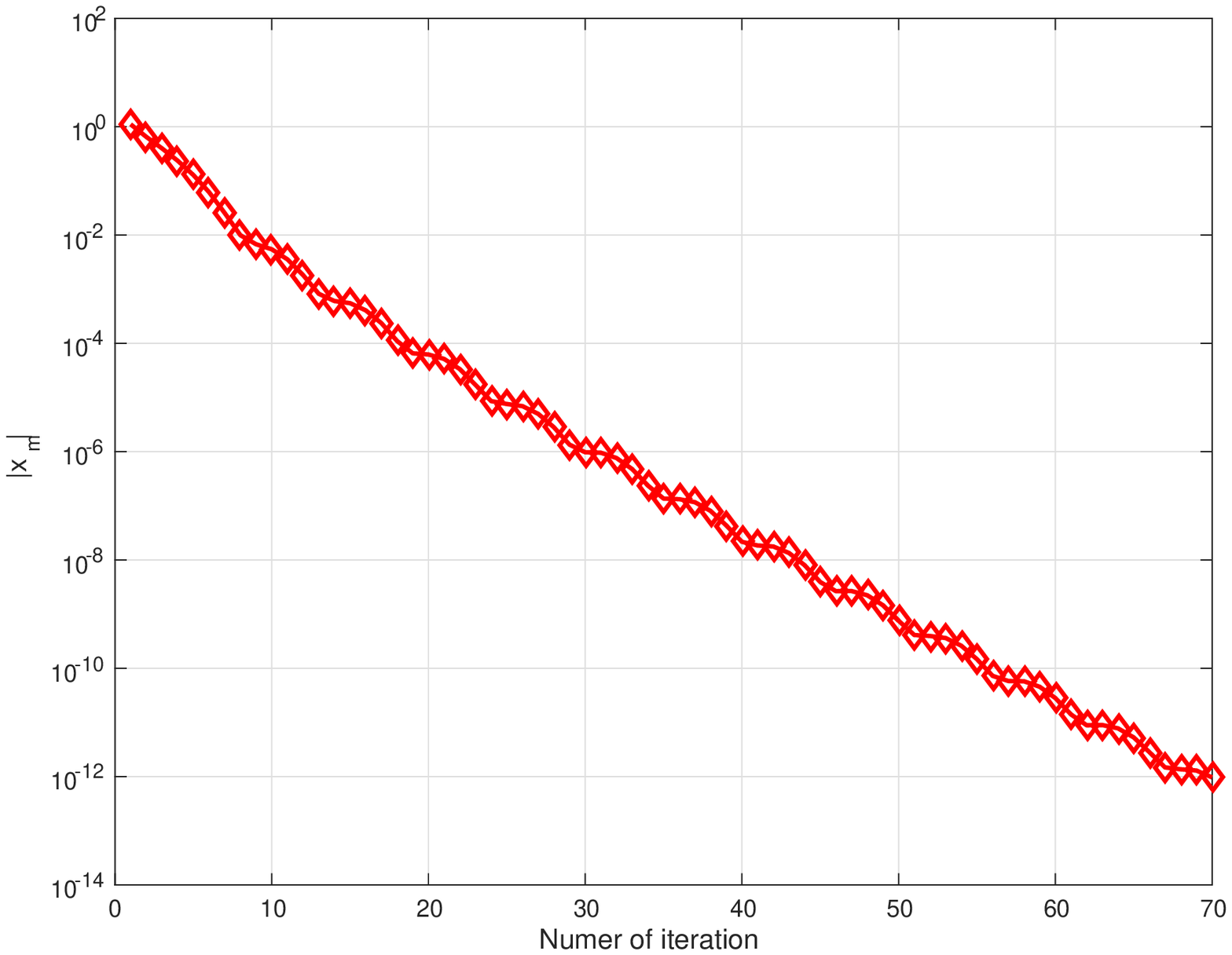}
 \includegraphics[width=0.460\textwidth]{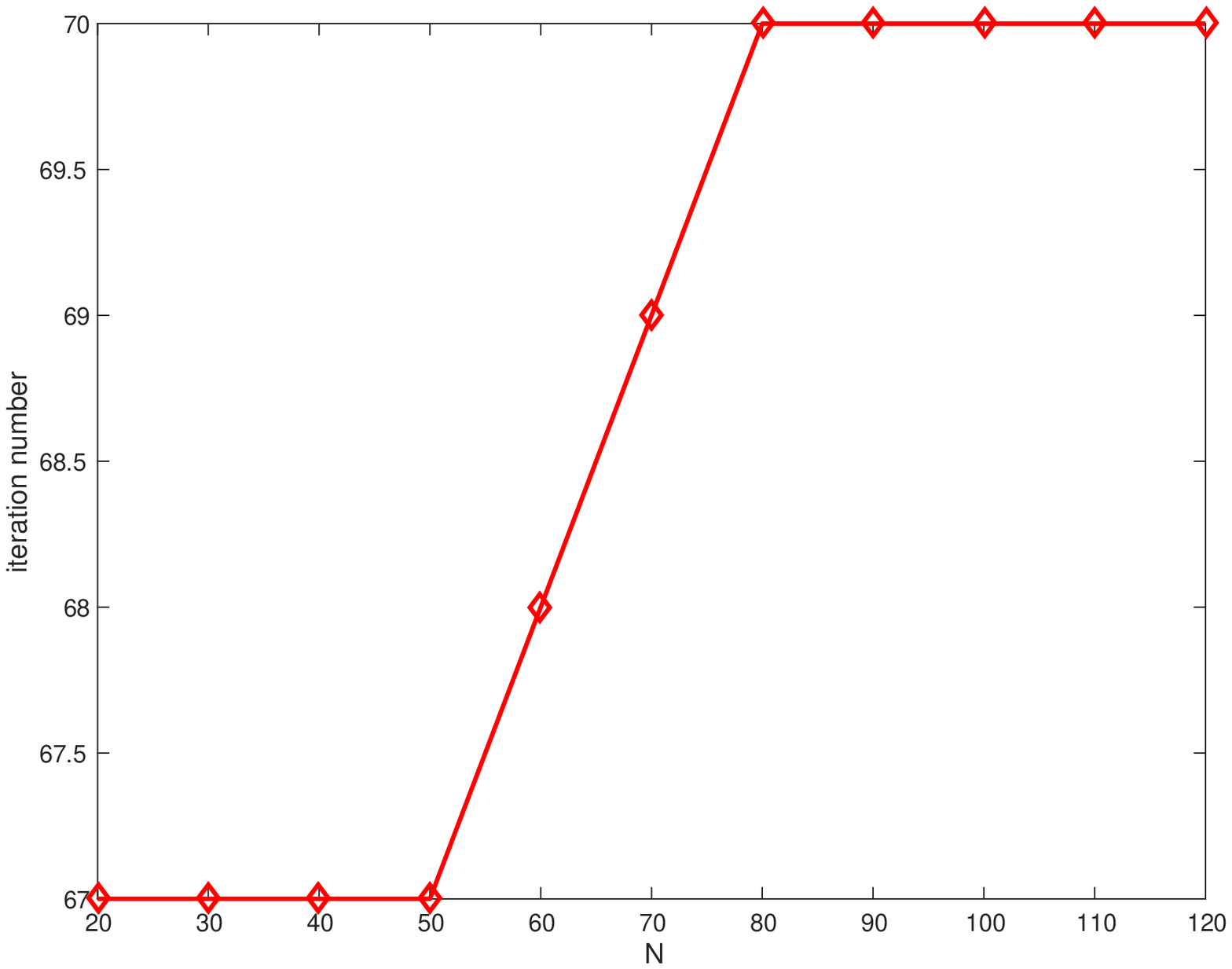}\\
 \centerline{Fig. 4.  Convergence of LaHOC iteration: Left) the error reduction in each iteration when $N=120$; } \\
 \centerline{Right) number of iterations needed when error reduction threshold is $10^{-12}$.}\\
 \end{tabular*}\\

  \noindent\textbf{ Test problem 3.2.}
  Consider the Euler dynamics and kinematics of a rigid body related to control laws to regulate the attitude of spacecraft and aircraft \cite{Tang1}:
  \begin{equation}
    \label{optimal attitude1}
    \begin{cases}
      \dot \rho(t)=\frac{1}{2}(I-S(\rho(t))+\rho(t)\rho^T(t))\omega(t), \\
      \dot \omega(t)=J^{-1}S(\omega(t))\,J\,\omega(t)+J^{-1}u(t),
    \end{cases}
  \end{equation}
  where $J=diag (10,6.3,8.5),$
  $\rho=(\rho_1,\rho_2,\rho_3)^T \in \mathbb{R}^3$ is the
  vector of Rodrigues parameters,
  $\omega=(\omega_1,\omega_2,\omega_3)^T \in \mathbb{R}^3,$
  is the angular velocity, and $u=(u_1,u_2,u_3)^T \in
  \mathbb{R}^3,$ is the control torque. The symbol $S(.)$ is
  a skew symmetric matrix of the form 
  \begin{equation}
    \label{S}
    S(\omega)  = \left[\begin{array}{ccc}
        0          & \omega_3 &  -\omega_2\\
        -\omega_3  & 0        &  \omega_1\\
        \omega_2  & -\omega_1 &  0
      \end{array}\right],
  \end{equation}
  In addition, the initial conditions are $\rho(0)=(0.3735, 0.4115, 0.2521)^T$ and $\omega(0)=(0,0,0)^T.$\\
  Then, according to the optimal control theory
  (\ref{Eq:Optimality conditions 2}), the optimality
  conditions can be written as:
  \begin{eqnarray}
    \label{system3}
    &&\dot \rho_1(t)=\frac{1}{2}\omega_1(t)+\frac{1}{2}\omega_1(t)\rho_1^2(t)+\frac{1}{2}\omega_2(t)\rho_1(t)\rho_2(t)+\frac{1}{2}\omega_3(t)\rho_1(t)\rho_3(t),\\
    &&\dot \rho_2(t)=\frac{1}{2}\omega_2(t)+\frac{1}{2}\omega_2(t)\rho_2^2(t)+\frac{1}{2}\omega_1(t)\rho_1(t)\rho_2(t)+\frac{1}{2}\omega_3(t)\rho_2(t)\rho_3(t),\\
    &&\dot \rho_3(t)=\frac{1}{2}\omega_3(t)+\frac{1}{2}\omega_3(t)\rho_3^2(t)+\frac{1}{2}\omega_1(t)\rho_1(t)\rho_3(t)+\frac{1}{2}\omega_2(t)\rho_2(t)\rho_3(t),\\
    &&\dot \omega_1(t)=-\frac{11}{50}\omega_2(t)\omega_3(t)-\frac{1}{100}\lambda_4(t),\\
    &&\dot \omega_2(t)=-\frac{5}{21}\omega_1(t)\omega_3(t)-\frac{100}{3969}\lambda_5(t),\\
    &&\dot \omega_3(t)=\frac{37}{85}\omega_1(t)\omega_2(t)-\frac{4}{289}\lambda_6(t),\\
    &&\dot \lambda_1(t)=-\lambda_1(t)\omega_1(t)\rho_1(t)-\rho_1(t)-\frac{1}{2}\lambda_1(t)\omega_2(t)\rho_2(t)-\frac{1}{2}\lambda_1(t)\omega_3(t)\rho_3(t)\nonumber\\
    &&~~~~~~~~~-\frac{1}{2}\lambda_2(t)\omega_1(t)\rho_2(t)-\frac{1}{2}\lambda_3(t)\omega_1(t)\rho_3(t),\\
    &&\dot \lambda_2(t)=-\lambda_2(t)\omega_2(t)\rho_2(t)-\rho_2(t)-\frac{1}{2}\lambda_1(t)\omega_2(t)\rho_1(t)-\frac{1}{2}\lambda_2(t)\omega_1(t)\rho_1(t)\nonumber\\
    &&~~~~~~~~~-\frac{1}{2}\lambda_2(t)\omega_3(t)\rho_3(t)-\frac{1}{2}\lambda_3(t)\omega_2(t)\rho_3(t),\\
    &&\dot \lambda_3(t)=-\lambda_3(t)\omega_3(t)\rho_3(t)-\rho_3(t)-\frac{1}{2}\lambda_1(t)\omega_3(t)\rho_1(t)-\frac{1}{2}\lambda_2(t)\omega_3(t)\rho_2(t)\nonumber\\
    &&~~~~~~~~~-\frac{1}{2}\lambda_3(t)\omega_1(t)\rho_1(t)-\frac{1}{2}\lambda_3(t)\omega_2(t)\rho_2(t),\\
    &&\dot \lambda_4(t)=-\frac{37}{85}\lambda_6(t)\omega_2(t)+\frac{5}{21}\lambda_5(t)\omega_3(t)-\frac{1}{2}\lambda_1(t)\rho_1^2(t)-\frac{1}{2}\lambda_2(t)\rho_1(t)\rho_2(t)\nonumber\\
    &&~~~~~~~~~-\frac{1}{2}\lambda_3(t)\rho_1(t)\rho_3(t)-\frac{1}{2}\lambda_1(t)-\omega_1(t),
      \end{eqnarray}
     \begin{eqnarray}
    &&\dot \lambda_5(t)=\frac{11}{50}\lambda_4(t)\omega_3(t)-\frac{37}{85}\lambda_6(t)\omega_1(t)-\frac{1}{2}\lambda_2(t)\rho_2^2(t)-\frac{1}{2}\lambda_1(t)\rho_1(t)\rho_2(t)\nonumber\\
    &&~~~~~~~~~-\frac{1}{2}\lambda_3(t)\rho_2(t)\rho_3(t)-\frac{1}{2}\lambda_2(t)-\omega_2(t),\\
    &&\dot \lambda_6(t)=-\frac{11}{50}\lambda_4(t)\omega_2(t)+\frac{5}{21}\lambda_5(t)\omega_1(t)-\frac{1}{2}\lambda_3(t)\rho_3^2(t)-\frac{1}{2}\lambda_1(t)\rho_1(t)\rho_3(t)\nonumber\\
    &&~~~~~~~~~-\frac{1}{2}\lambda_2(t)\rho_2(t)\rho_3(t)-\frac{1}{2}\lambda_3(t)-\omega_3(t),\\
    &&\rho_1(0)=0.3735, \quad \rho_2(0)= 0.4115,\quad \rho_3(0)= 0.2521,\quad \omega_1(0)=0, \quad \omega_2(0)=0, \quad\omega_3(0)=0,\nonumber
  \end{eqnarray}
  and the optimal control laws are $u_1(t)=-\frac{1}{10}\lambda_4,~ u_2(t)=-\frac{10}{63}\lambda_5,~ u_3(t)=-\frac{2}{17}\lambda_6.$\\
  In this example, the parameters used in the LaHOC
  algorithms are

  \begin{align}
    &\mathcal{L}_r = \left[\begin{array}{cccccccccccc}
        \frac{d}{dt} & 0 & 0 & -\frac{1}{2} & 0 & 0 & 0 & 0 & 0 & 0 & 0 & 0\\
        0 & \frac{d}{dt} & 0 & 0 & -\frac{1}{2} & 0 & 0 & 0 & 0 & 0 & 0 & 0\\
        0   &  0 & \frac{d}{dt} & 0 & 0 & -\frac{1}{2} & 0 & 0 & 0 & 0 & 0 & 0\\
        0   &  0 & 0 &  \frac{d}{dt} & 0 & 0 & 0 & 0 & 0 & \frac{1}{100} & 0 & 0\\
        0   &  0 & 0 & 0 & \frac{d}{dt} & 0 & 0 & 0 & 0 & 0 & \frac{100}{3969} & 0\\
        0   &  0 & 0 &  0 & 0 & \frac{d}{dt} &0 & 0 & 0 & 0 & 0 & \frac{4}{289}\\
        1   &  0 & 0 & 0 & 0 & 0 & \frac{d}{dt} & 0 & 0 & 0 & 0 & 0\\
        0   &  1 & 0 & 0 & 0 & 0 & 0 & \frac{d}{dt} & 0 & 0 & 0 & 0\\
        0   &  0 & 1 & 0 & 0 & 0 & 0 & 0 &  \frac{d}{dt} & 0 & 0 & 0\\
        0   &  0 & 0 & 1 & 0 & 0 & \frac{1}{2} & 0 & 0 & \frac{d}{dt} & 0 & 0\\
        0   &  0 & 0 & 0 & 1 & 0 & 0 & \frac{1}{2} & 0 & 0 & \frac{d}{dt} & 0\\
        0   &  0 & 0 & 0 & 0 & 1 & 0 & 0 & \frac{1}{2} & 0 & 0 & \frac{d}{dt}\\
      \end{array}\right],\\
       &\textbf{A} = \left[\begin{array}{cccccccccccc}
        \textbf{D}& \textbf{O} & \textbf{O}  & -\frac{1}{2}I & \textbf{O} & \textbf{O}  & \textbf{O} & \textbf{O} & \textbf{O}  & \textbf{O} & \textbf{O} & \textbf{O}  \\
        \textbf{O} &  \textbf{D} & \textbf{O} & \textbf{O} & -\frac{1}{2}I & I  & \textbf{O} & \textbf{O} & \textbf{O}  & \textbf{O} & \textbf{O} & \textbf{O}  \\
        \textbf{O} & \textbf{O} &  \textbf{D} & \textbf{O}& \textbf{O} & -\frac{1}{2}I  & \textbf{O} & \textbf{O} & \textbf{O}  & \textbf{O} & \textbf{O} & \textbf{O} \\
        \textbf{O} & \textbf{O} &  \textbf{O} & \textbf{D} & \textbf{O} & \textbf{O}   & \textbf{O} & \textbf{O} & \textbf{O}  & \frac{1}{100}I & \textbf{O} & \textbf{O}  \\
        \textbf{O} & \textbf{O} &  \textbf{O} & \textbf{O} & \textbf{D}  & \textbf{O}  & \textbf{O} & \textbf{O} & \textbf{O}  & \textbf{O} & \frac{100}{3969}I & \textbf{O}  \\
        \textbf{O} & \textbf{O} &  \textbf{O} & \textbf{O} & \textbf{O}  & \textbf{D}  & \textbf{O} & \textbf{O}  & \textbf{O} & \textbf{O} & \textbf{O} & \frac{4}{289}I \\
        I               & \textbf{O} &  \textbf{O} & \textbf{O} & \textbf{O}  & \textbf{O}  & \textbf{D} & \textbf{O} & \textbf{O}  & \textbf{O} & \textbf{O} & \textbf{O}  \\
        \textbf{O} & I               &  \textbf{O} & \textbf{O} & \textbf{O}  & \textbf{O}  & \textbf{O} & \textbf{D} & \textbf{O}  & \textbf{O} & \textbf{O} & \textbf{O}  \\
        \textbf{O} & \textbf{O} &  I               & \textbf{O} & \textbf{O}  & \textbf{O} & \textbf{O} & \textbf{O}  & \textbf{D} &\textbf{O} & \textbf{O} & \textbf{O}  \\
        \textbf{O} & \textbf{O} &  \textbf{O} &  I & \textbf{O}  & \textbf{O} & \frac{1}{2}I & \textbf{O}  & \textbf{O} & \textbf{D} & \textbf{O} & \textbf{O}  \\
        \textbf{O} & \textbf{O} &  \textbf{O} & \textbf{O} & I  & \textbf{O} & \textbf{O} & \frac{1}{2}I  & \textbf{O} & \textbf{O} &  \textbf{D} & \textbf{O}  \\
        \textbf{O} & \textbf{O} &  \textbf{O} &  \textbf{O} & \textbf{O}  & I & \textbf{O} & \textbf{O} & \frac{1}{2}I & \textbf{O} & \textbf{O}  &\textbf{D}\\
      \end{array}\right],
      \end{align}
    \begin{align}
    &R_{r,m-1} = \mathcal{L}_r[x_{r,m-1}] + Q_{r,m-1},\\
    &Q_{r,m-1} = \left[\begin{array}{l}
        \frac{1}{2}\displaystyle\sum_{j=0}^{m-1}\textbf{Z}_{4,m-1-j}\displaystyle\sum_{k=0}^{j}\textbf{Z}_{1,j}\textbf{Z}_{1,j-k}+\displaystyle\sum_{j=0}^{m-1}\textbf{Z}_{5,m-1-j}\displaystyle\sum_{k=0}^{j}\textbf{Z}_{1,j}\textbf{Z}_{2,j-k}\\
        +\displaystyle\sum_{j=0}^{m-1}\textbf{Z}_{6,m-1-j}\displaystyle\sum_{k=0}^{j}\textbf{Z}_{1,j}\textbf{Z}_{3,j-k},\\
        \frac{1}{2}\displaystyle\sum_{j=0}^{m-1}\textbf{Z}_{5,m-1-j}\displaystyle\sum_{k=0}^{j}\textbf{Z}_{2,j}\textbf{Z}_{2,j-k}+\displaystyle\sum_{j=0}^{m-1}\textbf{Z}_{4,m-1-j}\displaystyle\sum_{k=0}^{j}\textbf{Z}_{1,j}\textbf{Z}_{2,j-k}\\
        +\displaystyle\sum_{j=0}^{m-1}\textbf{Z}_{6,m-1-j}\displaystyle\sum_{k=0}^{j}\textbf{Z}_{2,j}\textbf{Z}_{3,j-k},\\
        \frac{1}{2}\displaystyle\sum_{j=0}^{m-1}\textbf{Z}_{6,m-1-j}\displaystyle\sum_{k=0}^{j}\textbf{Z}_{3,j}\textbf{Z}_{3,j-k}+\displaystyle\sum_{j=0}^{m-1}\textbf{Z}_{4,m-1-j}\displaystyle\sum_{k=0}^{j}\textbf{Z}_{1,j}\textbf{Z}_{3,j-k}\\
        +\displaystyle\sum_{j=0}^{m-1}\textbf{Z}_{5,m-1-j}\displaystyle\sum_{k=0}^{j}\textbf{Z}_{2,j}\textbf{Z}_{3,j-k},\\
      \end{array}\right],\quad r=1,2,3,
  \end{align}
  \begin{align}
  Q_{r,m-1} =  \left[\begin{array}{l}             
      -\frac{11}{50}\displaystyle \sum_{j=0}^{m-1}\textbf{Z}_{5,j} \textbf{Z}_{6,m-1-j},\\
      -\frac{5}{21}\displaystyle \sum_{j=0}^{m-1}\textbf{Z}_{4,j} \textbf{Z}_{6,m-1-j},\\
      \frac{37}{85}\displaystyle \sum_{j=0}^{m-1}\textbf{Z}_{4,j} \textbf{Z}_{5,m-1-j},\\
      -\displaystyle\sum_{j=0}^{m-1}\textbf{Z}_{7,m-1-j}\displaystyle\sum_{k=0}^{j}\textbf{Z}_{4,j}\textbf{Z}_{1,j-k}-\frac{1}{2}\displaystyle\sum_{j=0}^{m-1}\textbf{Z}_{7,m-1-j}\displaystyle\sum_{k=0}^{j}\textbf{Z}_{5,j}\textbf{Z}_{2,j-k}-\frac{1}{2}\displaystyle\sum_{j=0}^{m-1}\textbf{Z}_{7,m-1-j}\\
      \displaystyle\sum_{k=0}^{j}\textbf{Z}_{6,j}\textbf{Z}_{3,j-k}-\frac{1}{2}\displaystyle\sum_{j=0}^{m-1}\textbf{Z}_{8,m-1-j}\displaystyle\sum_{k=0}^{j}\textbf{Z}_{4,j}\textbf{Z}_{2,j-k}-\frac{1}{2}\displaystyle\sum_{j=0}^{m-1}\textbf{Z}_{9,m-1-j}\displaystyle\sum_{k=0}^{j}\textbf{Z}_{4,j}\textbf{Z}_{3,j-k},\\
      -\displaystyle\sum_{j=0}^{m-1}\textbf{Z}_{8,m-1-j}\displaystyle\sum_{k=0}^{j}\textbf{Z}_{5,j}\textbf{Z}_{2,j-k}-\frac{1}{2}\displaystyle\sum_{j=0}^{m-1}\textbf{Z}_{7,m-1-j}\displaystyle\sum_{k=0}^{j}\textbf{Z}_{5,j}\textbf{Z}_{1,j-k}-\frac{1}{2}\displaystyle\sum_{j=0}^{m-1}\textbf{Z}_{8,m-1-j}\\
      \displaystyle\sum_{k=0}^{j}\textbf{Z}_{4,j}\textbf{Z}_{1,j-k}-\frac{1}{2}\displaystyle\sum_{j=0}^{m-1}\textbf{Z}_{8,m-1-j}\displaystyle\sum_{k=0}^{j}\textbf{Z}_{6,j}\textbf{Z}_{3,j-k}-\frac{1}{2}\displaystyle\sum_{j=0}^{m-1}\textbf{Z}_{8,m-1-j}\displaystyle\sum_{k=0}^{j}\textbf{Z}_{5,j}\textbf{Z}_{3,j-k},\\
      -\displaystyle\sum_{j=0}^{m-1}\textbf{Z}_{9,m-1-j}\displaystyle\sum_{k=0}^{j}\textbf{Z}_{6,j}\textbf{Z}_{3,j-k}-\frac{1}{2}\displaystyle\sum_{j=0}^{m-1}\textbf{Z}_{7,m-1-j}\displaystyle\sum_{k=0}^{j}\textbf{Z}_{6,j}\textbf{Z}_{1,j-k}-\frac{1}{2}\displaystyle\sum_{j=0}^{m-1}\textbf{Z}_{8,m-1-j}\\
      \displaystyle\sum_{k=0}^{j}\textbf{Z}_{6,j}\textbf{Z}_{2,j-k}-\frac{1}{2}\displaystyle\sum_{j=0}^{m-1}\textbf{Z}_{9,m-1-j}\displaystyle\sum_{k=0}^{j}\textbf{Z}_{4,j}\textbf{Z}_{1,j-k}-\frac{1}{2}\displaystyle\sum_{j=0}^{m-1}\textbf{Z}_{9,m-1-j}\displaystyle\sum_{k=0}^{j}\textbf{Z}_{5,j}\textbf{Z}_{2,j-k},\\
      -\frac{37}{85}\displaystyle \sum_{j=0}^{m-1}\textbf{Z}_{12,j} \textbf{Z}_{5,m-1-j}+\frac{5}{21}\displaystyle \sum_{j=0}^{m-1}\textbf{Z}_{11,j} \textbf{Z}_{6,m-1-j}
      -\frac{1}{2}\displaystyle\sum_{j=0}^{m-1}\textbf{Z}_{7,m-1-j}\displaystyle\sum_{k=0}^{j}\textbf{Z}_{1,j}\textbf{Z}_{1,j-k}\\
      -\frac{1}{2}\displaystyle\sum_{j=0}^{m-1}\textbf{Z}_{8,m-1-j}\displaystyle\sum_{k=0}^{j}\textbf{Z}_{1,j}\textbf{Z}_{2,j-k}-\frac{1}{2}\displaystyle\sum_{j=0}^{m-1}\textbf{Z}_{9,m-1-j}\displaystyle\sum_{k=0}^{j}\textbf{Z}_{1,j}\textbf{Z}_{3,j-k},\\
      \frac{11}{50}\displaystyle \sum_{j=0}^{m-1}\textbf{Z}_{10,j} \textbf{Z}_{6,m-1-j}-\frac{37}{85}\displaystyle \sum_{j=0}^{m-1}\textbf{Z}_{12,j} \textbf{Z}_{4,m-1-j}
      -\frac{1}{2}\displaystyle\sum_{j=0}^{m-1}\textbf{Z}_{8,m-1-j}\displaystyle\sum_{k=0}^{j}\textbf{Z}_{2,j}\textbf{Z}_{2,j-k}\\
      -\frac{1}{2}\displaystyle\sum_{j=0}^{m-1}\textbf{Z}_{7,m-1-j}\displaystyle\sum_{k=0}^{j}\textbf{Z}_{1,j}\textbf{Z}_{2,j-k}-\frac{1}{2}\displaystyle\sum_{j=0}^{m-1}\textbf{Z}_{9,m-1-j}\displaystyle\sum_{k=0}^{j}\textbf{Z}_{2,j}\textbf{Z}_{3,j-k},\\
      -\frac{11}{50}\displaystyle \sum_{j=0}^{m-1}\textbf{Z}_{10,j} \textbf{Z}_{5,m-1-j}+\frac{5}{21}\displaystyle \sum_{j=0}^{m-1}\textbf{Z}_{11,j} \textbf{Z}_{4,m-1-j}-\frac{1}{2}\displaystyle\sum_{j=0}^{m-1}\textbf{Z}_{9,m-1-j}\displaystyle\sum_{k=0}^{j}\textbf{Z}_{3,j}\textbf{Z}_{3,j-k}\\
      -\frac{1}{2}\displaystyle\sum_{j=0}^{m-1}\textbf{Z}_{7,m-1-j}\displaystyle\sum_{k=0}^{j}\textbf{Z}_{1,j}\textbf{Z}_{3,j-k}
      -\frac{1}{2}\displaystyle\sum_{j=0}^{m-1}\textbf{Z}_{8,m-1-j}\displaystyle\sum_{k=0}^{j}\textbf{Z}_{2,j}\textbf{Z}_{3,j-k},\\
    \end{array}\right]\\
\end{align}
With these definitions, the LaHOC algorithm gives
\begin{equation}\label{soln1b}
  \textbf{X}_{r,m} =  (\chi_m + \hbar_r )\textbf{X}_{r,m-1}  + \hbar_r \textbf{A}^{-1}\textbf{Q}_{r,m-1},
\end{equation}
Because the right hand side of equation (\ref{soln1b}) is known, the solution can easily be obtained by using methods for solving linear system of equations.

Tables 2 and 3, give a comparison between the present LaHOC results for  $N = 50$ and $\hbar = -1$ and the numerically generated $\mathtt{BVP5C}$  at selected values of time $t.$ It can be seen from the tables that there is in good agreement between the two results.
\att{Moreover, our calculations show that the accuracy of LaHOC is faster}. In comparison with the $\mathtt{BVP5C}$, it is noteworthy that the LaHOC controls the error bounds while preserving the CPU time. The CPU time of LaHOC is $1.009860~s$, and $\mathtt{BVP5C}$ is $4.514071~s.$

Figurs. 5-9 show the suboptimal states and control for $m = 20$  iterations
of LaHOC, compared to  MATLAB built-in function $\mathtt{BVP5C}$. The convergence of Laguerre-LaHOC ieteration is depicted in Figure 10. 

The obtained optimal trajectories and optimal controls are almost identical to those obtained by
Jajarmi et al. \cite{Jajarmi}.


\begin{table}[!h]
 \caption{Comparison between the LaHOC solution when $N = 50$ and $\hbar = -1$  and $\mathtt{BVP5C}$  solution.}
       \label{Label}\centering
      \begin{tabular*}{0.99\textwidth}{@{\extracolsep{\fill}}ccccccc}\\\hline
     &  & $\rho_1(t)$  &  & $\rho_2(t)$ & & $\rho_3(t)$  \\
     \cline{2-3}\cline{4-5}\cline{6-7}
       \\
  \hline   
  $t$    & $LaHOC$ & $BVP5C$ & $LaHOC$ & $BVP5C$ & $LaHOC$ & $BVP5C$ \\
        \hline
0.409 	 & 0.371513 	 & 0.371389 	 & 0.408328 	 & 0.408146 	 & 0.251619 	 & 0.250403 \\ 
1.950 	 & 0.337343 	 & 0.335574 	 & 0.355885 	 & 0.353367 	 & 0.241531 	 & 0.221942 \\ 
4.663 	 & 0.237026 	 & 0.232989 	 & 0.215281 	 & 0.210198 	 & 0.198296 	 & 0.144524 \\ 
8.597 	 & 0.107445 	 & 0.103268 	 & 0.066722 	 & 0.062265 	 & 0.112940 	 & 0.057554 \\ 
20.488 	 & -0.010891 	 & -0.011225 	 & -0.006986 	 & -0.007030 	 & -0.003736 	 & -0.005378 \\ 
38.855 	 & 0.000248 	 & 0.000274 	 & 0.000140 	 & 0.000138 	 & 0.000053 	 & 0.000159 \\ \hline        
\end{tabular*}
       \end{table}

\begin{table}[!h]
 \caption{Comparison between the LaHOC solution when $N = 50$ and $\hbar = -1$  and $\mathtt{BVP5C}$  solution.}
       \label{Label}\centering
      \begin{tabular*}{0.99\textwidth}{@{\extracolsep{\fill}}ccccccc}\\\hline
     &  & $\omega_1(t)$  &  & $\omega_2(t)$ & & $\omega_3(t)$  \\
     \cline{2-3}\cline{4-5}\cline{6-7}
       \\\midrule
  \hline   
  $t$    & $LaHOC$ & $BVP5C$ & $LaHOC$ & $BVP5C$ & $LaHOC$ & $BVP5C$ \\
\hline
0.409 	 & -0.013313 	 & -0.013421 	 & -0.023872 	 & -0.024420 	 & -0.000899 	 & -0.011641 \\
1.950 	 & -0.047399 	 & -0.047871 	 & -0.077873 	 & -0.079134 	 & -0.009730 	 & -0.039806 \\
4.663 	 & -0.066195 	 & -0.067190 	 & -0.090208 	 & -0.090958 	 & -0.032073 	 & -0.049396 \\
8.597 	 & -0.051563 	 & -0.051458 	 & -0.050159 	 & -0.049573 	 & -0.040398 	 & -0.031829 \\
20.488 	 & -0.000271 	 & 0.000040 	 & 0.002386 	 & 0.002574 	 & -0.002575 	 & 0.000386 \\
38.855 	 & 0.000147 	 & 0.000132 	 & -0.000059 	 & -0.000062 	 & 0.000127 	 & 0.000009 \\
\hline \end{tabular*}
       \end{table}

\noindent\begin{tabular*}{.5cm}{cc}
  \includegraphics[scale=0.45]{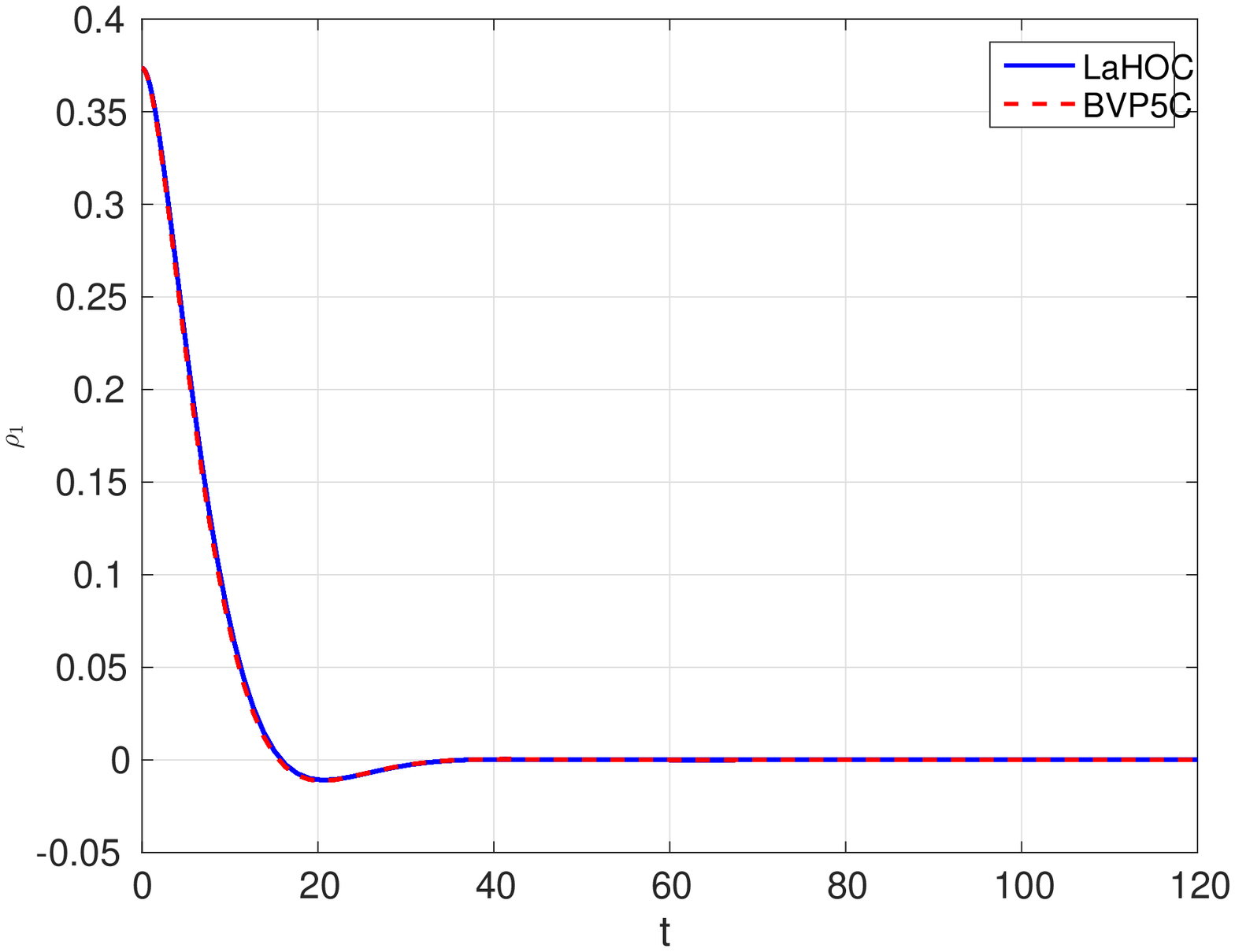} \includegraphics[scale=0.45]{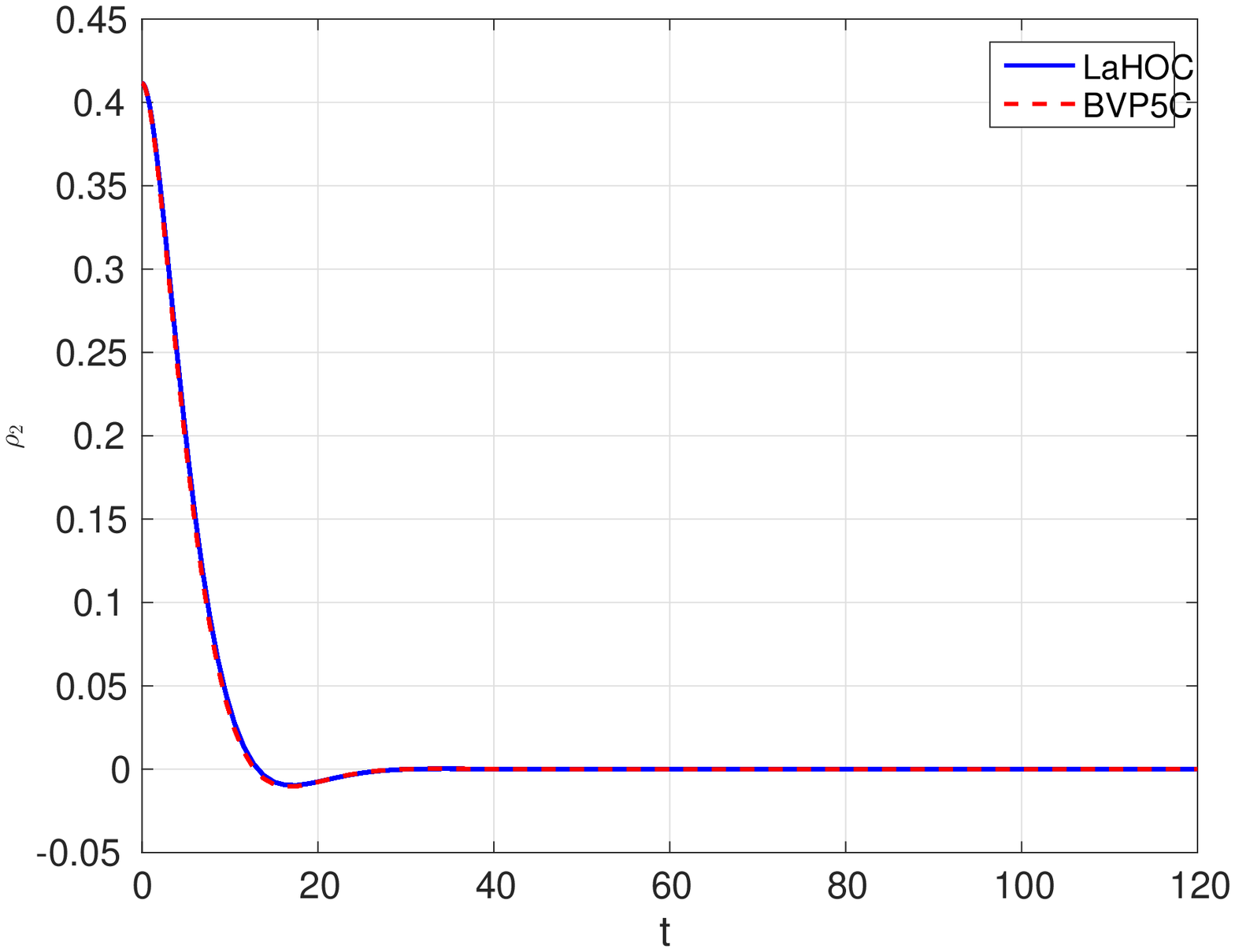}\\
  \centerline{Fig. 5. The amplitudes of optimal state variables.}  \\
\end{tabular*}\\
\begin{tabular*}{.5cm}{cc}
  \includegraphics[scale=0.45]{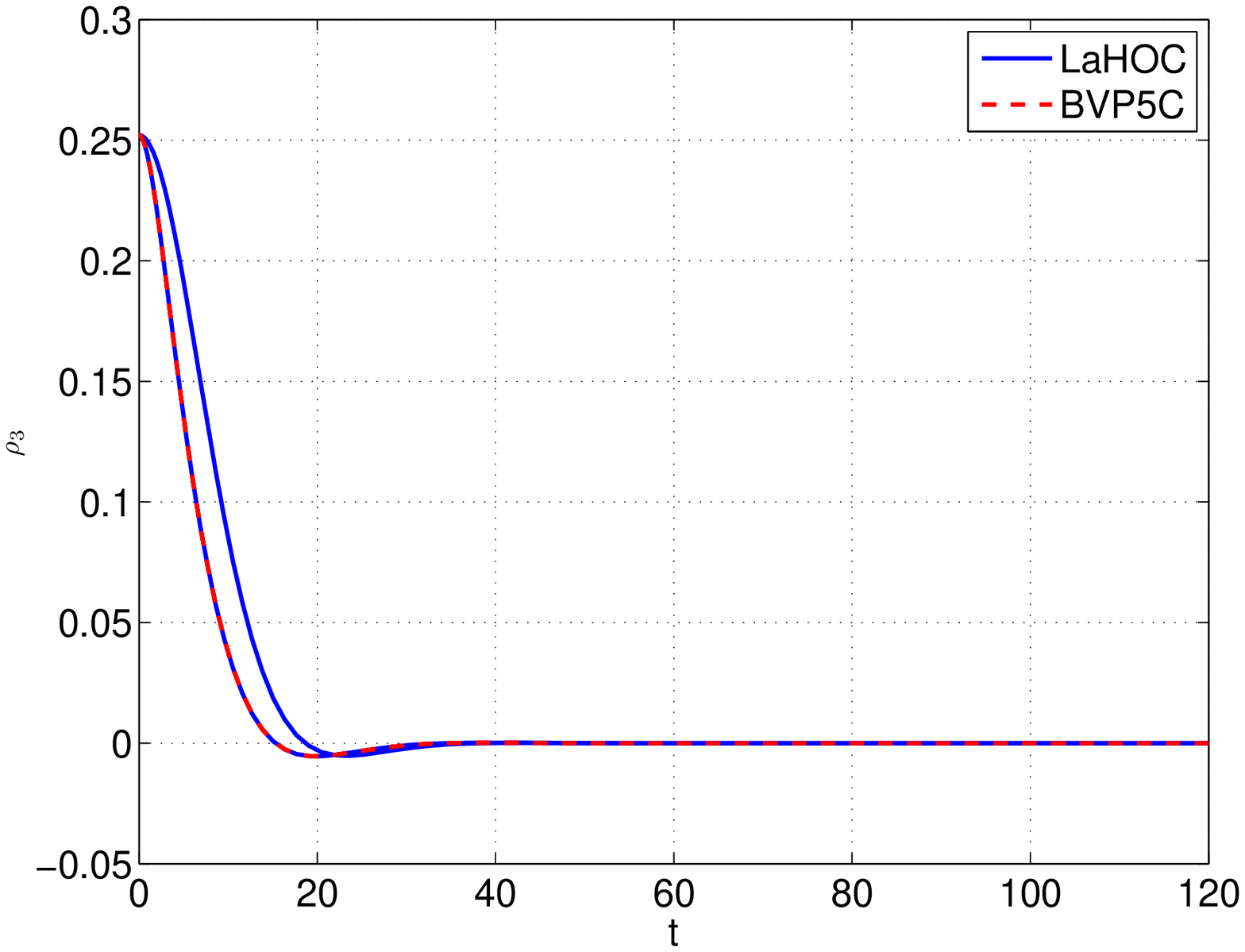} \includegraphics[scale=0.45]{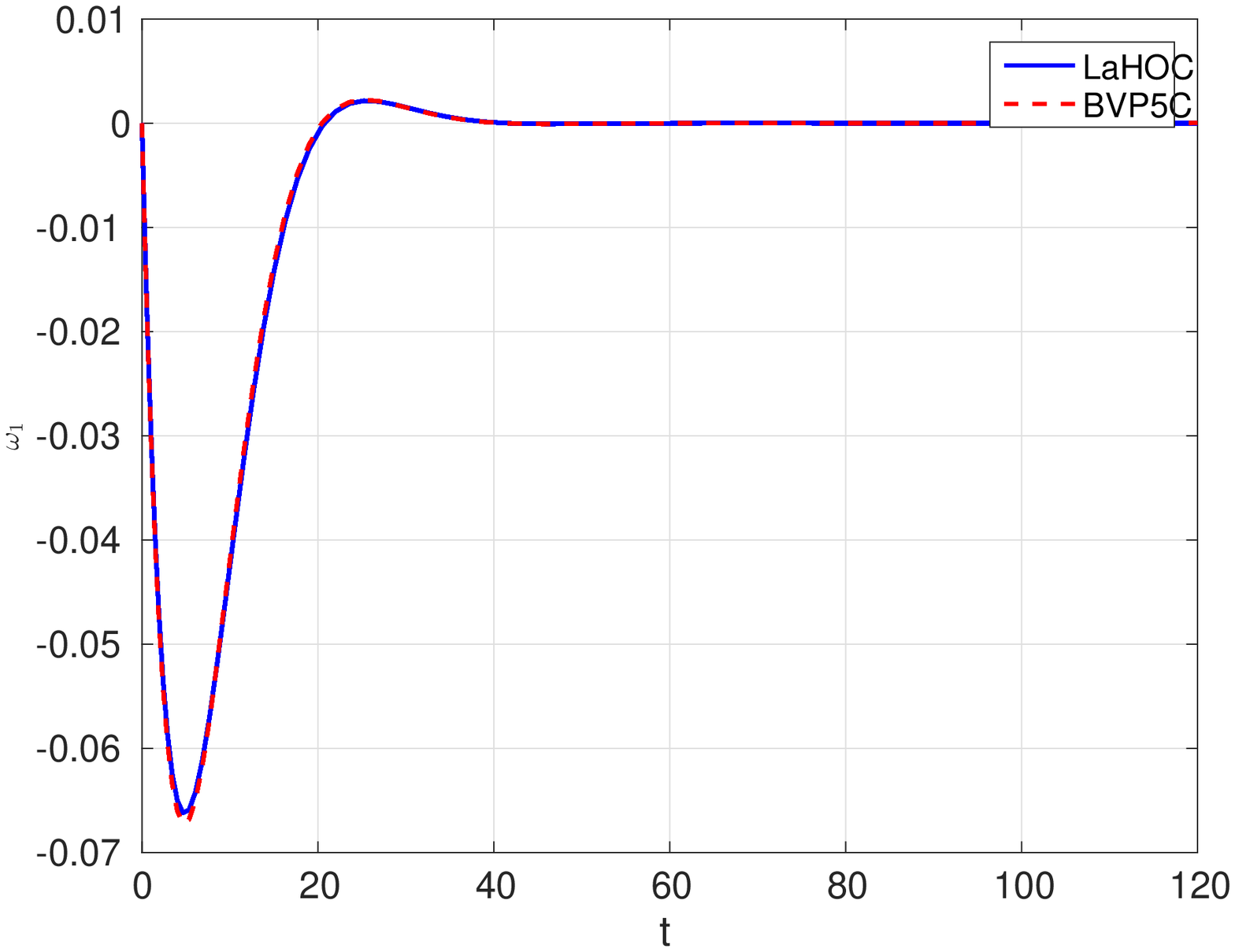}\\
  \centerline{Fig. 6. The amplitudes of optimal state variables.}  \\
\end{tabular*}\\
\begin{tabular*}{.5cm}{cc}
  \includegraphics[scale=0.45]{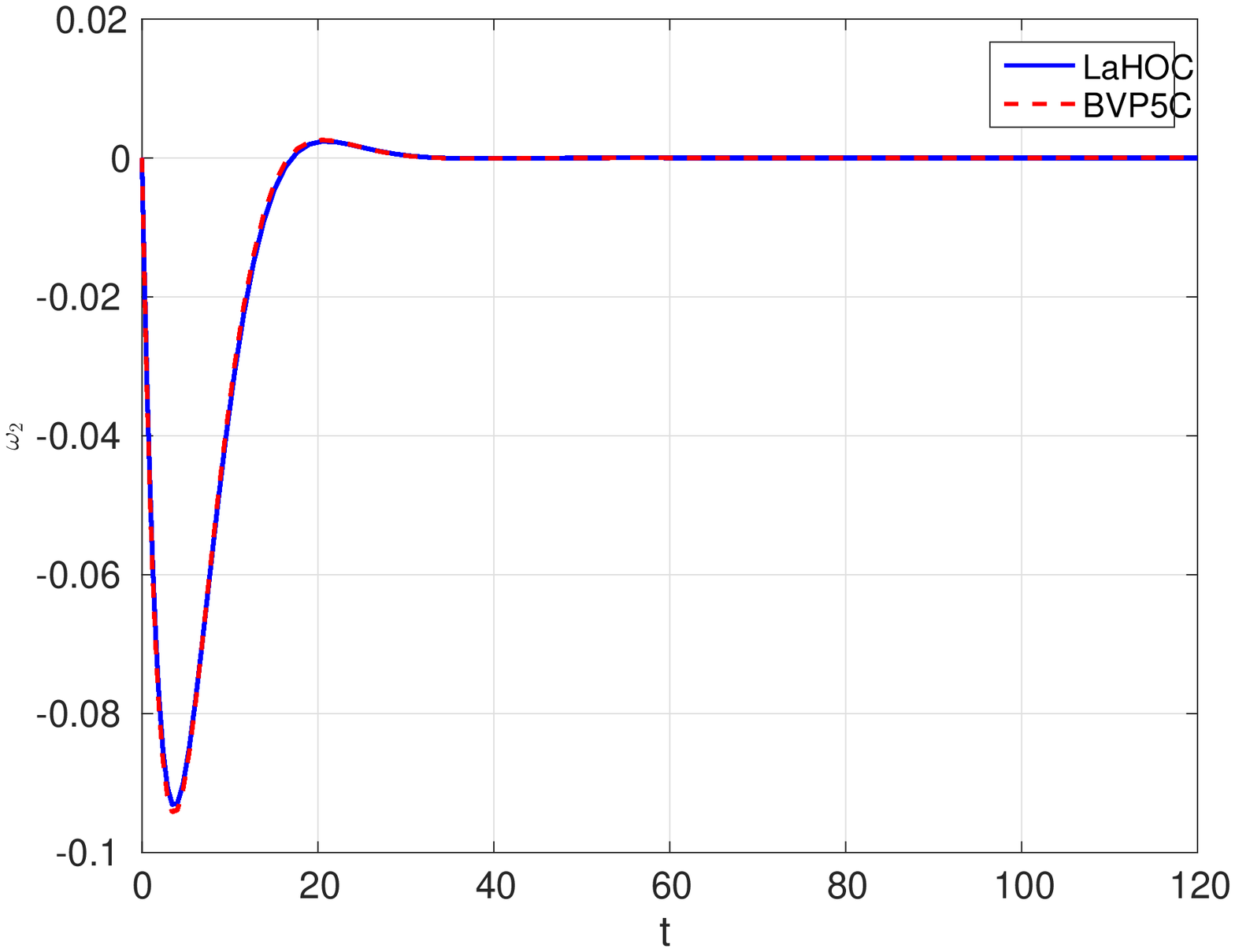} \includegraphics[scale=0.45]{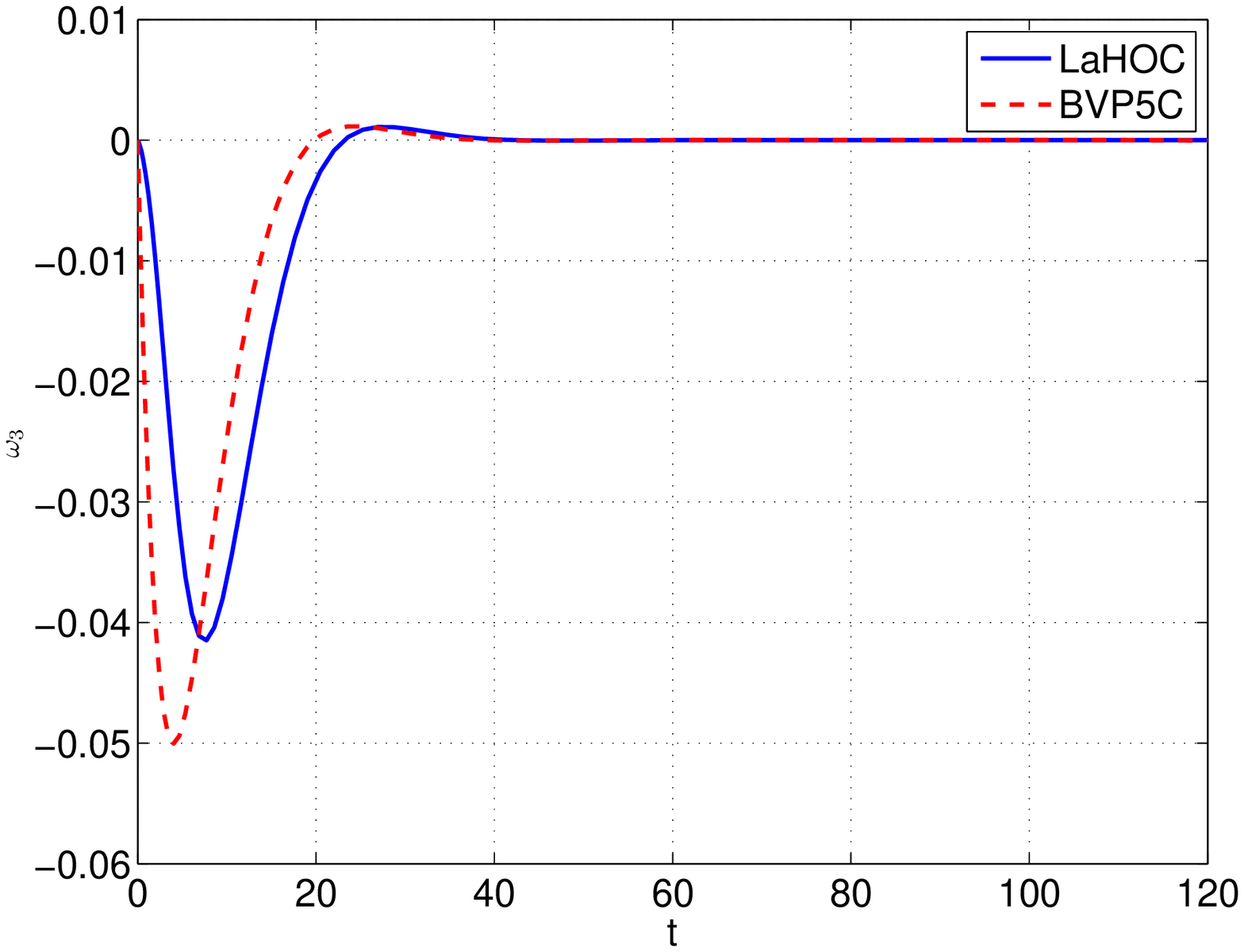}\\
  \centerline{Fig. 7. The amplitudes of optimal state variables.}  \\
\end{tabular*}\\

\noindent
\begin{tabular*}{.5cm}{cc}
  \includegraphics[scale=0.45]{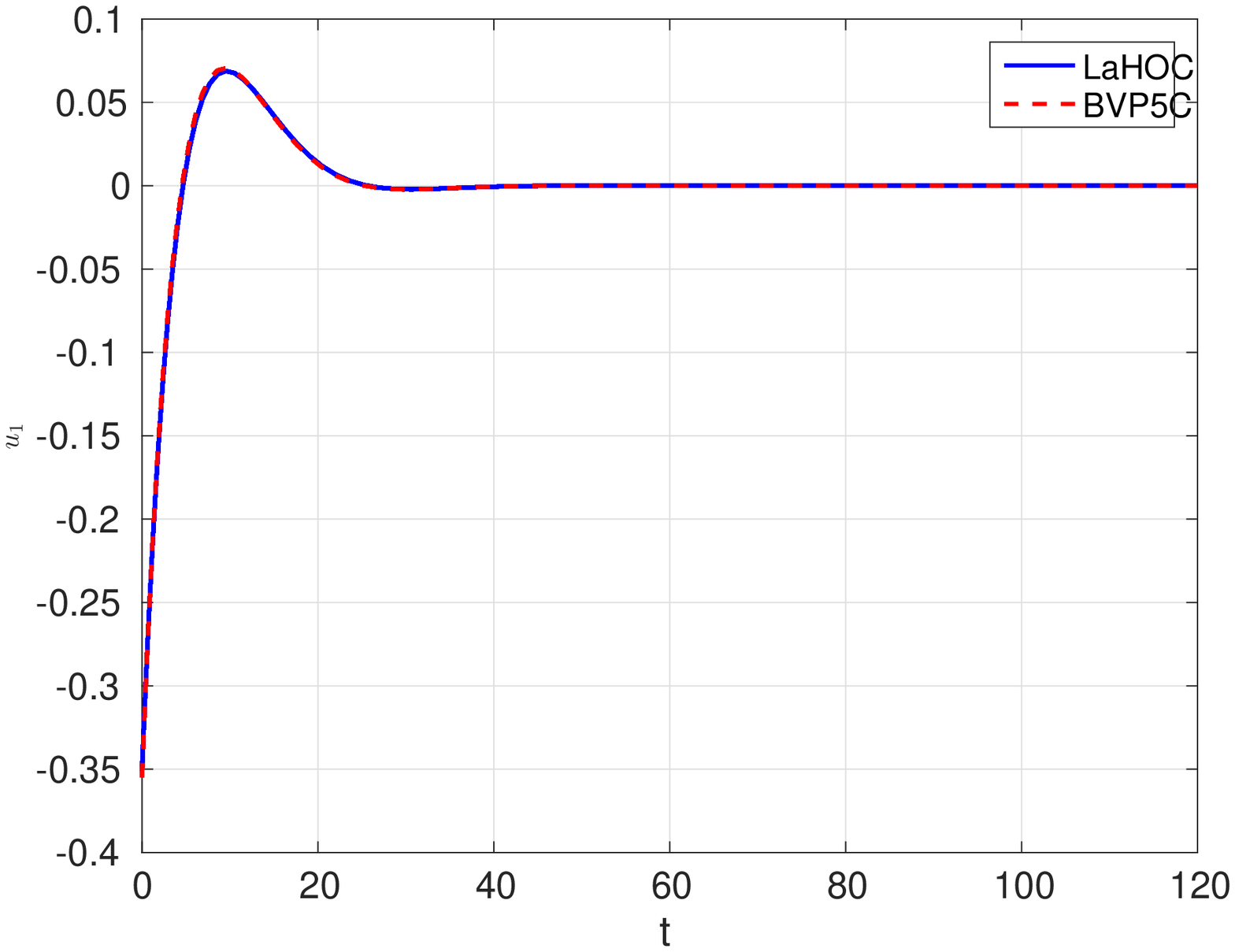} \includegraphics[scale=0.45]{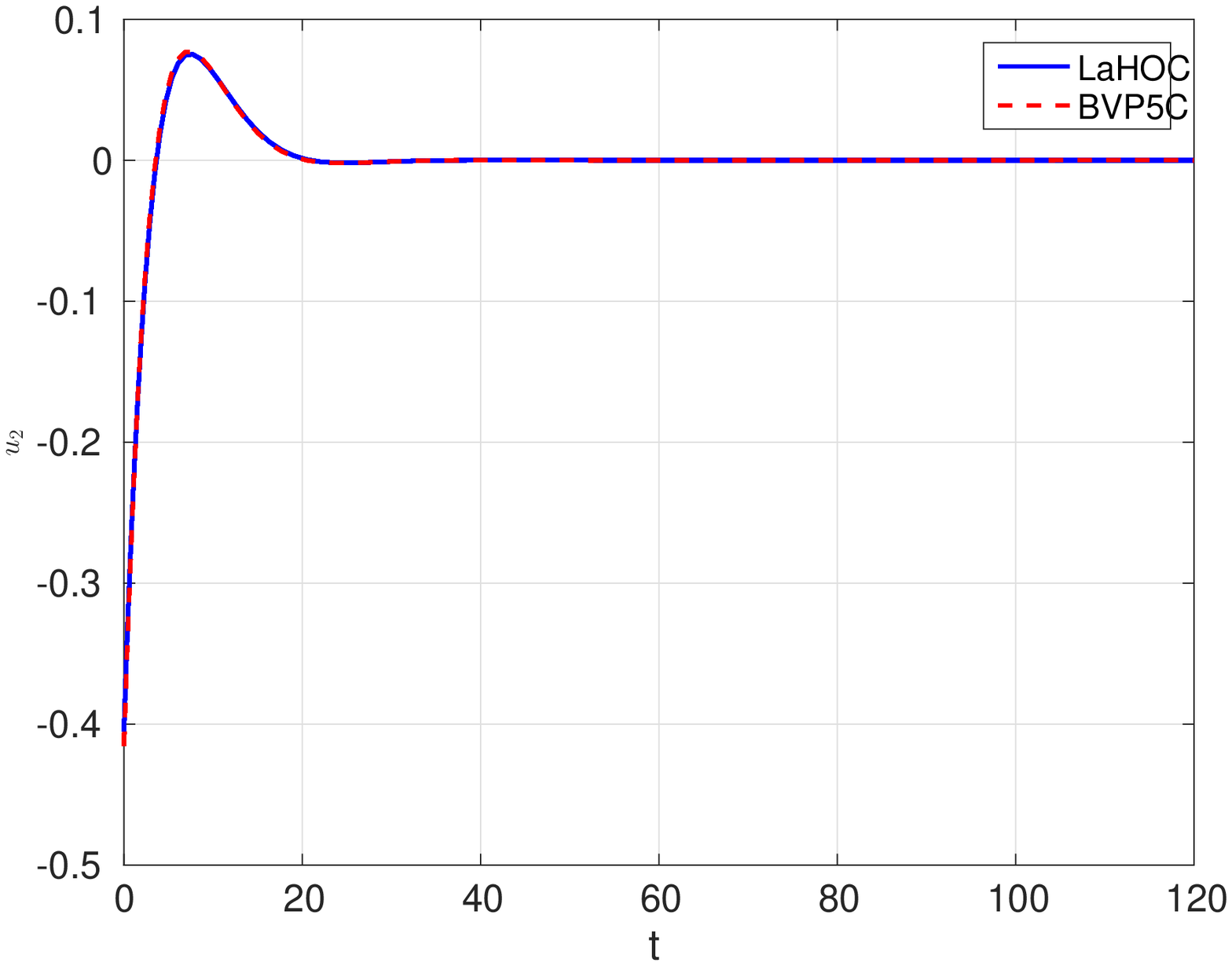}\\
  \centerline{Fig. 8. The amplitudes of optimal control variables.}  \\
\end{tabular*}\\
\begin{tabular*}{.5cm}{cc}
  \includegraphics[scale=0.45]{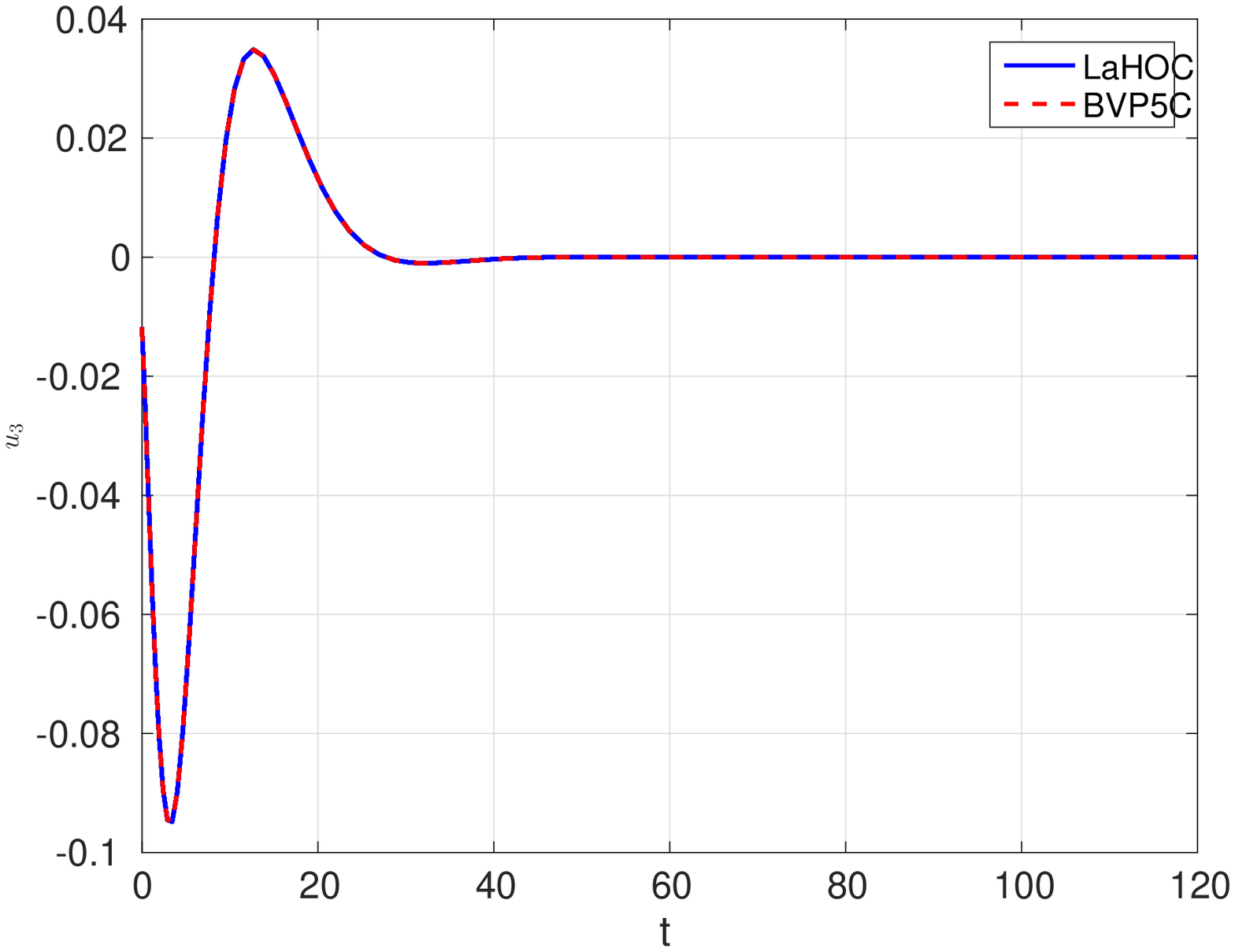}\\
  \centerline{Fig. 9. The amplitudes of optimal control variables.}  \\
\end{tabular*}\\
\begin{tabular*}{.5cm}{cc}
  \includegraphics[scale=0.45]{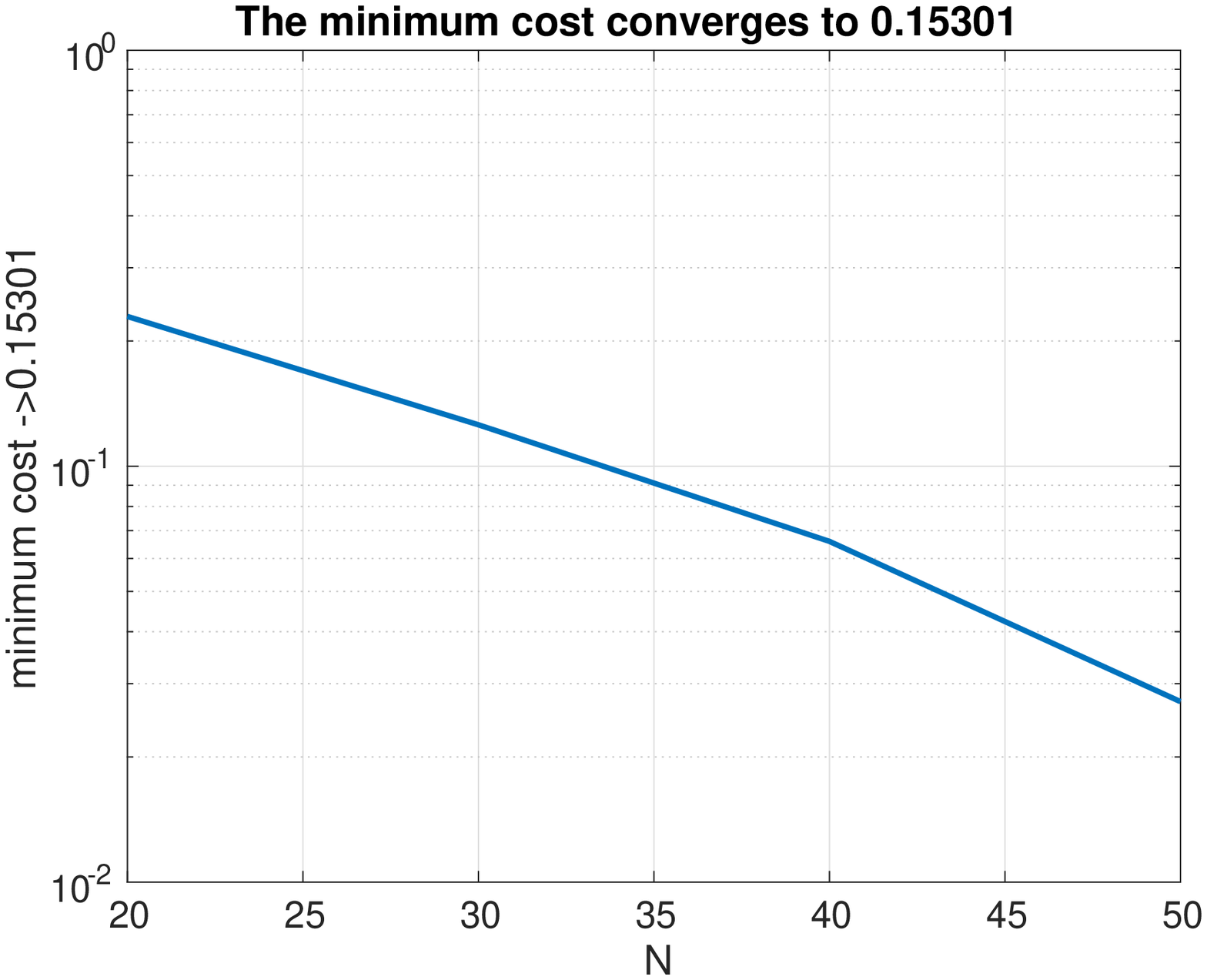}\\
 \centerline{Fig. 10.  The minimum cost onvergence.}  \\
 \end{tabular*}\\

\section{Conclusion}
In this paper, an effective method based upon the spectral homotopy method with Laguerre basis (LaHOC) is proposed for finding the numerical solutions of the  infinite horizon optimal control problem of nonlinear interconnected large-scale dynamic systems. Modified Laguerre method is used to discretize the equation of optimal condition, while homotopy method is used to construct an iterative scheme.
Two illustrative examples demonstrated that LaHOC has  spectral accuracy and  very good efficiency, which is comparable to well established numerical methods such as the MATLAB $\textsf{BVP5C}$ solver. The second example shows when the multi-components have different time and amplitude scales, one need to use 
adaptive rescaling technique in the Laguerre bases to improve accuracy, which deserves a further study.

%
%

\end{document}